\documentclass[12pt]{article}
\usepackage{amsfonts,amssymb}
\usepackage{url}
\usepackage{graphicx,color}
\usepackage{amsmath}
\newtheorem{defin}{Definition}[section]

\newtheorem{lem}[defin]{Lemma}
\newtheorem{lemma}[defin]{Lemma}

\newtheorem{thm}[defin]{Theorem}

\newtheorem{cor}[defin]{Corollary}

\def\qed{\hbox{\kern1pt\vrule height6pt width4pt
depth1pt\kern1pt}\medskip}
\def\bproof{\par\noindent{\bf Proof.\enspace}\rm}

\newcommand{\real}{{\mathbb{R}}}
\newcommand{\R}{{\mathbb{R}}}

\newcommand{\rat}{{\mathbb{Q}}}
\newcommand{\Z}{{\mathbb{Z}}}

\newcommand{\scrf}{{\mathcal{F}}}

\newcommand{\scri}{{\mathcal{I}}}

\newcommand{\scrm}{{\mathcal{M}}}

\newcommand{\scrp}{{\mathcal{P}}}
\newcommand{\scrq}{{\mathcal{Q}}}
\newcommand{\scrr}{{\mathcal{R}}}

\newcommand{\scrmcad}{{\mathcal{M}}_{\rm PL}}
\newcommand{\scrmcov}{{\mathcal{M}}_{\rm PL}}

\newcommand{\F}{{\cal F}}
\newcommand{\val}{{\rm val}}
\newcommand{\diag}{{\rm diag}}

\newcommand{\pl}{{PL}}

\newcommand{\sm}{\setminus}

\newcommand{\eproof}{\hfill $\bullet$\\}
\newcommand{\rank}{\mbox{\rm rank }}

\newcommand{\Null}{\mbox{\rm Null }}

\title{A characterisation of the generic rigidity of 2-dimensional
point-line frameworks}
\date{}
\author{Bill Jackson \thanks{School of Mathematical Sciences, Queen Mary
University of London, Mile End Road, London E1 4NS, United Kingdom. E-mail:
b.jackson@qmul.ac.uk} \and J.C. Owen \thanks{Siemens, Francis House,
112 Hills Road, Cambridge CB2 1PH, United Kingdom. E-mail:
owen.john.ext@siemens.com}}
\begin{document}
\maketitle

\begin{abstract}
A 2-dimensional point-line framework is a collection of points and
lines in the plane which are linked by pairwise constraints that fix
some angles between pairs of lines and also some point-line and
point-point distances. It is rigid if every continuous motion of the
points and lines which preserves the constraints results in a
point-line framework which can be obtained from the initial
framework by a translation or a rotation. We characterise when a
generic point-line framework is rigid. Our characterisation gives
rise to a polynomial algorithm for solving this decision problem.
\end{abstract}


\section{Introduction}

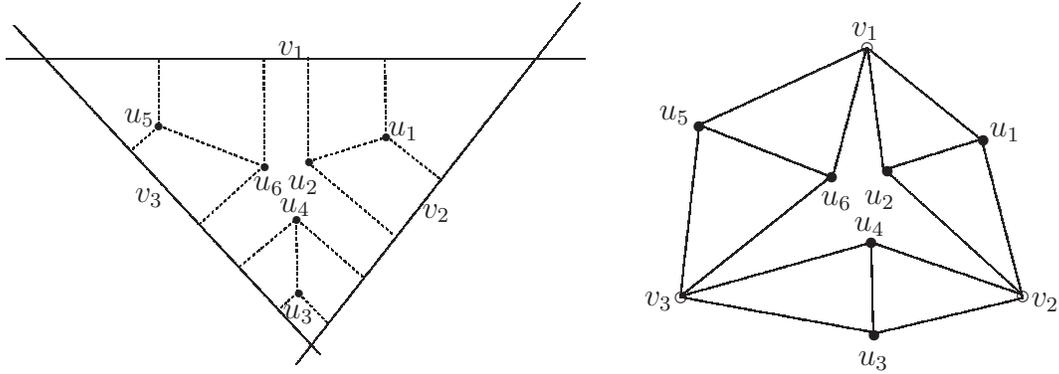
\begin{figure}[ht]
\begin{center}
\unitlength .4mm 
\linethickness{0.4pt}
\ifx\plotpoint\undefined\newsavebox{\plotpoint}\fi 
\begin{picture}(124.5,121.5)(50,-20)
\put(92.75,80.5){\circle*{3.041}}
\put(128.021,67.021){\circle*{3.041}}
\put(138.521,49.521){\circle*{3.041}}
\put(142.771,68.521){\circle*{3.041}}
\put(139.271,24.771){\circle*{3.041}}
\put(168.271,76.771){\circle*{3.041}}
\put(137,105.75){\makebox(0,0)[cc]{$v_1$}}
\put(89.5,56.75){\makebox(0,0)[cc]{$v_3$}}
\put(174,79){\makebox(0,0)[cc]{$u_1$}}
\put(140.75,61){\makebox(0,0)[cc]{$u_2$}}
\put(139,18){\makebox(0,0)[cc]{$u_3$}}
\put(138,53.75){\makebox(0,0)[cc]{$u_4$}}
\put(86.25,83.5){\makebox(0,0)[cc]{$u_5$}}
\put(129.25,60.75){\makebox(0,0)[cc]{$u_6$}}
\put(42,103){\line(1,0){192.5}}
\multiput(232.5,121.5)(-.06746599857,-.08607730852){1397}{\line(0,-1){.08607730852}}
\multiput(44.5,114.25)(.06746031746,-.07258597884){1512}{\line(0,-1){.07258597884}}
\put(92.859,80.859){\line(0,1){.9565}}
\put(92.859,82.772){\line(0,1){.9565}}
\put(92.859,84.686){\line(0,1){.9565}}
\put(92.859,86.599){\line(0,1){.9565}}
\put(92.859,88.512){\line(0,1){.9565}}
\put(92.859,90.425){\line(0,1){.9565}}
\put(92.859,92.338){\line(0,1){.9565}}
\put(92.859,94.251){\line(0,1){.9565}}
\put(92.859,96.164){\line(0,1){.9565}}
\put(92.859,98.077){\line(0,1){.9565}}
\put(92.859,99.99){\line(0,1){.9565}}
\put(92.859,101.903){\line(0,1){.9565}}
\multiput(92.609,80.359)(-.073413,-.063492){9}{\line(-1,0){.073413}}
\multiput(91.288,79.217)(-.073413,-.063492){9}{\line(-1,0){.073413}}
\multiput(89.967,78.074)(-.073413,-.063492){9}{\line(-1,0){.073413}}
\multiput(88.645,76.931)(-.073413,-.063492){9}{\line(-1,0){.073413}}
\multiput(87.324,75.788)(-.073413,-.063492){9}{\line(-1,0){.073413}}
\multiput(86.002,74.645)(-.073413,-.063492){9}{\line(-1,0){.073413}}
\multiput(84.681,73.502)(-.073413,-.063492){9}{\line(-1,0){.073413}}
\multiput(92.359,80.859)(.151709,-.058761){6}{\line(1,0){.151709}}
\multiput(94.18,80.154)(.151709,-.058761){6}{\line(1,0){.151709}}
\multiput(96,79.449)(.151709,-.058761){6}{\line(1,0){.151709}}
\multiput(97.821,78.744)(.151709,-.058761){6}{\line(1,0){.151709}}
\multiput(99.641,78.039)(.151709,-.058761){6}{\line(1,0){.151709}}
\multiput(101.462,77.334)(.151709,-.058761){6}{\line(1,0){.151709}}
\multiput(103.282,76.629)(.151709,-.058761){6}{\line(1,0){.151709}}
\multiput(105.103,75.924)(.151709,-.058761){6}{\line(1,0){.151709}}
\multiput(106.924,75.218)(.151709,-.058761){6}{\line(1,0){.151709}}
\multiput(108.744,74.513)(.151709,-.058761){6}{\line(1,0){.151709}}
\multiput(110.565,73.808)(.151709,-.058761){6}{\line(1,0){.151709}}
\multiput(112.385,73.103)(.151709,-.058761){6}{\line(1,0){.151709}}
\multiput(114.206,72.398)(.151709,-.058761){6}{\line(1,0){.151709}}
\multiput(116.026,71.693)(.151709,-.058761){6}{\line(1,0){.151709}}
\multiput(117.847,70.988)(.151709,-.058761){6}{\line(1,0){.151709}}
\multiput(119.667,70.282)(.151709,-.058761){6}{\line(1,0){.151709}}
\multiput(121.488,69.577)(.151709,-.058761){6}{\line(1,0){.151709}}
\multiput(123.308,68.872)(.151709,-.058761){6}{\line(1,0){.151709}}
\multiput(125.129,68.167)(.151709,-.058761){6}{\line(1,0){.151709}}
\multiput(126.949,67.462)(.151709,-.058761){6}{\line(1,0){.151709}}
\put(127.859,67.109){\line(0,1){.973}}
\put(127.846,69.055){\line(0,1){.973}}
\put(127.832,71.001){\line(0,1){.973}}
\put(127.819,72.947){\line(0,1){.973}}
\put(127.805,74.893){\line(0,1){.973}}
\put(127.792,76.839){\line(0,1){.973}}
\put(127.778,78.785){\line(0,1){.973}}
\put(127.765,80.731){\line(0,1){.973}}
\put(127.751,82.677){\line(0,1){.973}}
\put(127.738,84.623){\line(0,1){.973}}
\put(127.724,86.569){\line(0,1){.973}}
\put(127.711,88.515){\line(0,1){.973}}
\put(127.697,90.461){\line(0,1){.973}}
\put(127.684,92.407){\line(0,1){.973}}
\put(127.67,94.353){\line(0,1){.973}}
\put(127.657,96.299){\line(0,1){.973}}
\put(127.643,98.245){\line(0,1){.973}}
\put(127.63,100.19){\line(0,1){.973}}
\put(127.616,102.136){\line(0,1){.973}}
\multiput(127.609,67.359)(-.069355,-.062903){10}{\line(-1,0){.069355}}
\multiput(126.222,66.101)(-.069355,-.062903){10}{\line(-1,0){.069355}}
\multiput(124.835,64.843)(-.069355,-.062903){10}{\line(-1,0){.069355}}
\multiput(123.448,63.585)(-.069355,-.062903){10}{\line(-1,0){.069355}}
\multiput(122.061,62.327)(-.069355,-.062903){10}{\line(-1,0){.069355}}
\multiput(120.674,61.069)(-.069355,-.062903){10}{\line(-1,0){.069355}}
\multiput(119.287,59.811)(-.069355,-.062903){10}{\line(-1,0){.069355}}
\multiput(117.9,58.553)(-.069355,-.062903){10}{\line(-1,0){.069355}}
\multiput(116.513,57.295)(-.069355,-.062903){10}{\line(-1,0){.069355}}
\multiput(115.126,56.037)(-.069355,-.062903){10}{\line(-1,0){.069355}}
\multiput(113.738,54.779)(-.069355,-.062903){10}{\line(-1,0){.069355}}
\multiput(112.351,53.521)(-.069355,-.062903){10}{\line(-1,0){.069355}}
\multiput(110.964,52.263)(-.069355,-.062903){10}{\line(-1,0){.069355}}
\multiput(109.577,51.005)(-.069355,-.062903){10}{\line(-1,0){.069355}}
\multiput(108.19,49.747)(-.069355,-.062903){10}{\line(-1,0){.069355}}
\multiput(106.803,48.488)(-.069355,-.062903){10}{\line(-1,0){.069355}}
\multiput(138.859,25.109)(-.069444,-.061111){9}{\line(-1,0){.069444}}
\multiput(137.609,24.009)(-.069444,-.061111){9}{\line(-1,0){.069444}}
\multiput(136.359,22.909)(-.069444,-.061111){9}{\line(-1,0){.069444}}
\multiput(135.109,21.809)(-.069444,-.061111){9}{\line(-1,0){.069444}}
\multiput(133.859,20.709)(-.069444,-.061111){9}{\line(-1,0){.069444}}
\multiput(139.109,25.109)(.065,-.066667){10}{\line(0,-1){.066667}}
\multiput(140.409,23.776)(.065,-.066667){10}{\line(0,-1){.066667}}
\multiput(141.709,22.443)(.065,-.066667){10}{\line(0,-1){.066667}}
\multiput(143.009,21.109)(.065,-.066667){10}{\line(0,-1){.066667}}
\multiput(144.309,19.776)(.065,-.066667){10}{\line(0,-1){.066667}}
\multiput(145.609,18.443)(.065,-.066667){10}{\line(0,-1){.066667}}
\multiput(146.909,17.109)(.065,-.066667){10}{\line(0,-1){.066667}}
\multiput(148.209,15.776)(.065,-.066667){10}{\line(0,-1){.066667}}
\put(138.359,49.359){\line(0,-1){.97}}
\put(138.399,47.419){\line(0,-1){.97}}
\put(138.439,45.479){\line(0,-1){.97}}
\put(138.479,43.539){\line(0,-1){.97}}
\put(138.519,41.599){\line(0,-1){.97}}
\put(138.559,39.659){\line(0,-1){.97}}
\put(138.599,37.719){\line(0,-1){.97}}
\put(138.639,35.779){\line(0,-1){.97}}
\put(138.679,33.839){\line(0,-1){.97}}
\put(138.719,31.899){\line(0,-1){.97}}
\put(138.759,29.959){\line(0,-1){.97}}
\put(138.799,28.019){\line(0,-1){.97}}
\put(138.839,26.079){\line(0,-1){.97}}
\multiput(138.109,49.359)(-.075,-.063){10}{\line(-1,0){.075}}
\multiput(136.609,48.099)(-.075,-.063){10}{\line(-1,0){.075}}
\multiput(135.109,46.839)(-.075,-.063){10}{\line(-1,0){.075}}
\multiput(133.609,45.579)(-.075,-.063){10}{\line(-1,0){.075}}
\multiput(132.109,44.319)(-.075,-.063){10}{\line(-1,0){.075}}
\multiput(130.609,43.059)(-.075,-.063){10}{\line(-1,0){.075}}
\multiput(129.109,41.799)(-.075,-.063){10}{\line(-1,0){.075}}
\multiput(127.609,40.539)(-.075,-.063){10}{\line(-1,0){.075}}
\multiput(126.109,39.279)(-.075,-.063){10}{\line(-1,0){.075}}
\multiput(124.609,38.019)(-.075,-.063){10}{\line(-1,0){.075}}
\multiput(123.109,36.759)(-.075,-.063){10}{\line(-1,0){.075}}
\multiput(121.609,35.499)(-.075,-.063){10}{\line(-1,0){.075}}
\multiput(120.109,34.239)(-.075,-.063){10}{\line(-1,0){.075}}
\multiput(138.609,49.609)(.070968,-.062097){10}{\line(1,0){.070968}}
\multiput(140.029,48.367)(.070968,-.062097){10}{\line(1,0){.070968}}
\multiput(141.448,47.126)(.070968,-.062097){10}{\line(1,0){.070968}}
\multiput(142.867,45.884)(.070968,-.062097){10}{\line(1,0){.070968}}
\multiput(144.287,44.642)(.070968,-.062097){10}{\line(1,0){.070968}}
\multiput(145.706,43.4)(.070968,-.062097){10}{\line(1,0){.070968}}
\multiput(147.126,42.158)(.070968,-.062097){10}{\line(1,0){.070968}}
\multiput(148.545,40.916)(.070968,-.062097){10}{\line(1,0){.070968}}
\multiput(149.964,39.674)(.070968,-.062097){10}{\line(1,0){.070968}}
\multiput(151.384,38.432)(.070968,-.062097){10}{\line(1,0){.070968}}
\multiput(152.803,37.19)(.070968,-.062097){10}{\line(1,0){.070968}}
\multiput(154.222,35.948)(.070968,-.062097){10}{\line(1,0){.070968}}
\multiput(155.642,34.706)(.070968,-.062097){10}{\line(1,0){.070968}}
\multiput(157.061,33.464)(.070968,-.062097){10}{\line(1,0){.070968}}
\multiput(158.48,32.222)(.070968,-.062097){10}{\line(1,0){.070968}}
\multiput(159.9,30.98)(.070968,-.062097){10}{\line(1,0){.070968}}
\put(142.359,68.609){\line(0,1){.9929}}
\put(142.359,70.595){\line(0,1){.9929}}
\put(142.359,72.581){\line(0,1){.9929}}
\put(142.359,74.567){\line(0,1){.9929}}
\put(142.359,76.552){\line(0,1){.9929}}
\put(142.359,78.538){\line(0,1){.9929}}
\put(142.359,80.524){\line(0,1){.9929}}
\put(142.359,82.509){\line(0,1){.9929}}
\put(142.359,84.495){\line(0,1){.9929}}
\put(142.359,86.481){\line(0,1){.9929}}
\put(142.359,88.467){\line(0,1){.9929}}
\put(142.359,90.452){\line(0,1){.9929}}
\put(142.359,92.438){\line(0,1){.9929}}
\put(142.359,94.424){\line(0,1){.9929}}
\put(142.359,96.409){\line(0,1){.9929}}
\put(142.359,98.395){\line(0,1){.9929}}
\put(142.359,100.381){\line(0,1){.9929}}
\put(142.359,102.367){\line(0,1){.9929}}
\multiput(141.859,68.859)(.073125,-.061875){10}{\line(1,0){.073125}}
\multiput(143.322,67.622)(.073125,-.061875){10}{\line(1,0){.073125}}
\multiput(144.784,66.384)(.073125,-.061875){10}{\line(1,0){.073125}}
\multiput(146.247,65.147)(.073125,-.061875){10}{\line(1,0){.073125}}
\multiput(147.709,63.909)(.073125,-.061875){10}{\line(1,0){.073125}}
\multiput(149.172,62.672)(.073125,-.061875){10}{\line(1,0){.073125}}
\multiput(150.634,61.434)(.073125,-.061875){10}{\line(1,0){.073125}}
\multiput(152.097,60.197)(.073125,-.061875){10}{\line(1,0){.073125}}
\multiput(153.559,58.959)(.073125,-.061875){10}{\line(1,0){.073125}}
\multiput(155.022,57.722)(.073125,-.061875){10}{\line(1,0){.073125}}
\multiput(156.484,56.484)(.073125,-.061875){10}{\line(1,0){.073125}}
\multiput(157.947,55.247)(.073125,-.061875){10}{\line(1,0){.073125}}
\multiput(159.409,54.009)(.073125,-.061875){10}{\line(1,0){.073125}}
\multiput(160.872,52.772)(.073125,-.061875){10}{\line(1,0){.073125}}
\multiput(162.334,51.534)(.073125,-.061875){10}{\line(1,0){.073125}}
\multiput(163.797,50.297)(.073125,-.061875){10}{\line(1,0){.073125}}
\multiput(165.259,49.059)(.073125,-.061875){10}{\line(1,0){.073125}}
\multiput(166.722,47.822)(.073125,-.061875){10}{\line(1,0){.073125}}
\multiput(168.184,46.584)(.073125,-.061875){10}{\line(1,0){.073125}}
\multiput(169.647,45.347)(.073125,-.061875){10}{\line(1,0){.073125}}
\put(168.109,77.109){\line(0,1){.963}}
\put(168.091,79.035){\line(0,1){.963}}
\put(168.072,80.961){\line(0,1){.963}}
\put(168.054,82.887){\line(0,1){.963}}
\put(168.035,84.813){\line(0,1){.963}}
\put(168.017,86.739){\line(0,1){.963}}
\put(167.998,88.665){\line(0,1){.963}}
\put(167.98,90.591){\line(0,1){.963}}
\put(167.961,92.517){\line(0,1){.963}}
\put(167.943,94.443){\line(0,1){.963}}
\put(167.924,96.369){\line(0,1){.963}}
\put(167.906,98.295){\line(0,1){.963}}
\put(167.887,100.221){\line(0,1){.963}}
\put(167.869,102.146){\line(0,1){.963}}
\multiput(168.109,77.359)(.08,-.063333){9}{\line(1,0){.08}}
\multiput(169.549,76.219)(.08,-.063333){9}{\line(1,0){.08}}
\multiput(170.989,75.079)(.08,-.063333){9}{\line(1,0){.08}}
\multiput(172.429,73.939)(.08,-.063333){9}{\line(1,0){.08}}
\multiput(173.869,72.799)(.08,-.063333){9}{\line(1,0){.08}}
\multiput(175.309,71.659)(.08,-.063333){9}{\line(1,0){.08}}
\multiput(176.749,70.519)(.08,-.063333){9}{\line(1,0){.08}}
\multiput(178.189,69.379)(.08,-.063333){9}{\line(1,0){.08}}
\multiput(179.629,68.239)(.08,-.063333){9}{\line(1,0){.08}}
\multiput(181.069,67.099)(.08,-.063333){9}{\line(1,0){.08}}
\multiput(182.509,65.959)(.08,-.063333){9}{\line(1,0){.08}}
\multiput(183.949,64.819)(.08,-.063333){9}{\line(1,0){.08}}
\multiput(185.389,63.679)(.08,-.063333){9}{\line(1,0){.08}}
\multiput(142.359,68.859)(.183929,.057143){5}{\line(1,0){.183929}}
\multiput(144.199,69.431)(.183929,.057143){5}{\line(1,0){.183929}}
\multiput(146.038,70.002)(.183929,.057143){5}{\line(1,0){.183929}}
\multiput(147.877,70.574)(.183929,.057143){5}{\line(1,0){.183929}}
\multiput(149.717,71.145)(.183929,.057143){5}{\line(1,0){.183929}}
\multiput(151.556,71.717)(.183929,.057143){5}{\line(1,0){.183929}}
\multiput(153.395,72.288)(.183929,.057143){5}{\line(1,0){.183929}}
\multiput(155.234,72.859)(.183929,.057143){5}{\line(1,0){.183929}}
\multiput(157.074,73.431)(.183929,.057143){5}{\line(1,0){.183929}}
\multiput(158.913,74.002)(.183929,.057143){5}{\line(1,0){.183929}}
\multiput(160.752,74.574)(.183929,.057143){5}{\line(1,0){.183929}}
\multiput(162.592,75.145)(.183929,.057143){5}{\line(1,0){.183929}}
\multiput(164.431,75.717)(.183929,.057143){5}{\line(1,0){.183929}}
\multiput(166.27,76.288)(.183929,.057143){5}{\line(1,0){.183929}}
\put(185,51.75){\makebox(0,0)[cc]{$v_2$}}
\end{picture}
\unitlength .5mm 
\linethickness{0.4pt}
\ifx\plotpoint\undefined\newsavebox{\plotpoint}\fi 
\begin{picture}(164.75,105.75)(20,0)
\put(137.5,101.25){\circle{3.162}}
\put(87.831,34.831){\circle{3.162}}
\put(178.831,35.081){\circle{3.162}}
\put(92.75,80.5){\circle*{3.041}}
\put(128.021,67.021){\circle*{3.041}}
\put(138.521,49.521){\circle*{3.041}}
\put(142.771,68.521){\circle*{3.041}}
\put(139.271,24.771){\circle*{3.041}}
\put(168.271,76.771){\circle*{3.041}}
\multiput(92.75,80.75)(.170673077,-.067307692){208}{\line(1,0){.170673077}}
\multiput(128.25,66.75)(.067164179,.257462687){134}{\line(0,1){.257462687}}
\multiput(137.25,101.25)(-.1463815789,-.0674342105){304}{\line(-1,0){.1463815789}}
\multiput(92.75,80.75)(-.06666667,-.60666667){75}{\line(0,-1){.60666667}}
\multiput(87.75,35.25)(.239386792,.067216981){212}{\line(1,0){.239386792}}
\multiput(138.5,49.5)(.0625,-3.03125){8}{\line(0,-1){3.03125}}
\multiput(139,25.25)(.281914894,.067375887){141}{\line(1,0){.281914894}}
\multiput(178.75,34.75)(-.067073171,.257621951){164}{\line(0,1){.257621951}}
\multiput(167.75,77)(-.211382114,-.067073171){123}{\line(-1,0){.211382114}}
\multiput(141.75,68.75)(-.06746032,.51587302){63}{\line(0,1){.51587302}}
\multiput(137.5,101.25)(.0835579515,-.0673854447){371}{\line(1,0){.0835579515}}
\multiput(142,69)(.0740365112,-.0674442191){493}{\line(1,0){.0740365112}}
\multiput(178.5,35.75)(-.194711538,.067307692){208}{\line(-1,0){.194711538}}
\multiput(128,67.25)(-.0847280335,-.0674686192){478}{\line(-1,0){.0847280335}}
\multiput(87.5,35)(.368794326,-.067375887){141}{\line(1,0){.368794326}}
\put(137,105.75){\makebox(0,0)[cc]{$v_1$}}
\put(184.75,33.75){\makebox(0,0)[cc]{$v_2$}}
\put(82,34){\makebox(0,0)[cc]{$v_3$}}
\put(174,79){\makebox(0,0)[cc]{$u_1$}}
\put(140.75,61){\makebox(0,0)[cc]{$u_2$}}
\put(139,18){\makebox(0,0)[cc]{$u_3$}}
\put(138,53.75){\makebox(0,0)[cc]{$u_4$}}
\put(86.25,83.5){\makebox(0,0)[cc]{$u_5$}}
\put(129.25,60.75){\makebox(0,0)[cc]{$u_6$}}
\end{picture}

\end{center}
\caption{A point-line framework in $\R^2$ and its associated
point-line graph.} \label{fig1}
\end{figure}

A point-line framework is a collection of points and lines in
$d$-dimensional Euclidean space which are linked by pairwise
constraints that fix the angles between some pairs of lines, the
distances between some pairs of points and the distances between
some pairs of points and lines.  The placing of the pairwise
constraints is represented by a point-line graph
where the vertices in the graph correspond to the points and lines,
and an edge in the graph corresponds to the existence of a pairwise
constraint. A point-line framework is obtained from a point-line
graph by assigning coordinates to the points and lines, see Figure
\ref{fig1}.

A point-line framework is rigid if every continuous motion of the
points and lines which preserves the constraints results in a
point-line framework which can be obtained from the initial
framework by an isometry of the whole space.

The
constraints of a point-line framework determine its rigidity matrix,
in which rows are indexed by the edges of its point-line graph $G$.
The rigidity matroid of the framework is the matroid on the edge set
of $G$ defined by the linear independence  of the rows of the
rigidity matrix. Since this linear independence will be the same for
all generic frameworks with the same point-line graph, the rigidity
matroid of a generic framework is completely determined by the
dimension and the point-line graph. We will denote the
$2$-dimensional point-line rigidity matroid of a point-line graph
$G$ by $\scrm_{PL}(G)$ (or simply $\scrm_{PL}$ where the graph is
implied).

Point-line frameworks with no lines correspond to the much studied
bar-joint frameworks.
Such frameworks provide a model for a variety of physical systems
such as bar and joint structures \cite{WW} (where points correspond
to universal joints and bars correspond to distance constraints) or
molecular structures \cite{GH} (where points correspond to atoms and
distance constraints correspond to bonds).

Laman \cite{L} obtained the following characterisation of
independence in the generic 2-dimensional bar-joint rigidity
matroid. Given a graph $G=(V,E)$, let $\nu:2^E\to \Z$ by taking
$\nu(S)$ to be the number of vertices incident to $S$ for all
$S\subseteq E$. Then $S$ is independent in the generic 2-dimensional
bar-joint rigidity matroid of $G$ if and only if $|S'|\leq
2\nu(S')-3$ for all $\emptyset\neq S'\subseteq S$. The analogous
condition that $|S'|\leq d\nu(S')-d(d+1)/2$ is a necessary condition
for independence in the $d$-dimensional bar-joint rigidity matroid
but it is not sufficient when $d\geq 3$. Characterising independence
in the $d$-dimensional bar-joint rigidity matroid is an important
open problem. We refer the reader to the survey article of Whiteley
\cite{WW} for more information on bar-joint frameworks.

The point-line graph on the right of Figure \ref{fig1} can be used
to illustrate the difference between point-line and bar-joint
frameworks. We can use Laman's theorem to deduce that every generic
realisation of this graph as a 2-dimensional bar-joint framework
(with all vertices as points and all edges as distance constraints)
is rigid. In contrast we will see below that no generic realisation
as a 2-dimensional point-line framework is rigid.



Point-line graphs have been used extensively in computer aided
design \cite{O,OO,BF} and computer aided geometry \cite{WO,GC}. The
rigidity of the corresponding point-line framework determines when a
geometric design is well-dimensioned and is useful in determining
the decomposition of a geometric design into rigid components.


Our main result characterises independence in the 2-dimensional
point-line rigidity matroid $\scrm_{PL}$.
Our characterisation uses two count matroids\footnote{The term {\em
count matroid} refers to a matroid on the edge set of a graph
$G=(V,E)$, in which the independence of a set $S\subseteq E$ is
determined by counting the number of vertices incident to each
subset of $S$, see the Appendix of \cite{Wchapter} for more
details.} $\scrm(\nu_L)$ and $\scrm(2\nu_P+\nu_L-2)$ which are
defined as follows.
Given a point-line graph $G=(V,E)$, let $\nu_P,\nu_L:2^E\to \Z$ by
taking $\nu_P(S)$ and $\nu_L(S)$ to be the number of point-vertices,
respectively line-vertices, incident to $S$ for all $S\subseteq E$.
Then $S$ is independent in $\scrm(\nu_L)$, respectively
$\scrm(2\nu_P+\nu_L-2)$, if and only if $|S'|\leq \nu_L(S')$,
respectively  $|S'|\leq 2\nu_P(S')+\nu_L(S')-2$, for all
$\emptyset\neq S'\subseteq S$. We will show that $\scrm_{PL}$ is a
Dilworth truncation of the matroid union of $\scrm(\nu_L)$ and
$\scrm(2\nu_P+\nu_L-2)$.
Our characterisation of $\scrm_{PL}$ leads immediately to a
polynomial time algorithm to determine a maximal set of independent
edges in $\scrm_{PL}$.


We can consider a line as a one dimensional affine subspace of the
usual two dimensional Euclidean space. An alternative approach, which we use here, is to
consider a line as an oriented hyperplane with codimension one. When
$d=2$, the essential difference in the second approach is that the
distance from a point to a line becomes a signed quantity which
passes linearly through zero as the point crosses the line. This
corresponds to the coordinatisation of a line which is commonly used
in computer aided design \cite{OO} and is used by Yang \cite{Y} in
his extension of Cayley-Menger determinants and distance geometry to
include lines when $d=2$.  It is straightforward to show that these
two formulations lead to the same generic rigidity matroid when
$d=2$. However for $d\geq 3$ a line and an oriented hyperplane have
different dimensions and hence give rise to different rigidity
matrices and matroids. It remains an open problem to characterise the
generic rigidity matroid in either case when $d\geq 3$.

We have so far considered only the infinitesimal rigidity of
point-line frameworks. It may also be of interest to consider global
rigidity.  (Generic global rigidity for the special case of
$2$-dimensional bar-joint frameworks was characterised in
\cite{C,JJ}.) Note that it may be significant whether we view a line
as a one dimensional subspace or as an oriented hyperplane for
global rigidity even when $d=2$.

The remainder of this paper is organised as follows. In Section 2 we
give formal definitions for a point-line graph $G$, a 2-dimensional
point-line framework $(G,p)$, its rigidity matrix $R(G,p)$ and the
generic rigidity matroid $\scrm_{PL}(G)$. We also obtain two
necessary conditions for a rigidity matrix to have linearly
independent rows by considering its null space.

In Section 3 we use these two necessary conditions to define a
matroid $\scrm_\sharp$ on the edge set of a point-line graph with
the property that every independent set of $\scrm_{PL}$ is
independent in $\scrm_\sharp$. The derivation of this necessary
condition for independence in $\scrm_{PL}$ is purely combinatorial
(and should
extend to $d>2$). 
We
investigate the properties of $\scrm_\sharp$ and show that it can be
described as a Dilworth truncation of the matroid union of
$\scrm(\nu_L)$ and $\scrm(2\nu_P+\nu_L-2)$.

We complete our characterisation of $\scrm_{PL}$ in Section 4 by
showing that every independent set of $\scrm_\sharp$ is also
independent in $\scrm_{PL}$. Our result generalises Laman's Theorem
(which corresponds to the case when there are no line-vertices). Our
proof technique differs from that of Laman, however. Laman's proof
is based on a recursive construction of the family of graphs whose
edge set is independent in the generic bar-joint rigidity matroid.
We have not been able to obtain a similar recursive construction for
point-line graphs. Instead we adapt an alternative proof technique
for Laman's Theorem due to Whiteley \cite{W}, and give a direct
construction for a point-line framework $(G,p)$ with the property
that $E$ is independent in $\scrm_\pl(G,p)$ whenever $E$ is
independent in $\scrm_\sharp$.


In order to make this construction as simple as possible, we
restrict our attention to `naturally bipartite' point-line graphs,
i.e. point-line graphs in which every edge is incident with both a
point-vertex and a line-vertex. We then complete the proof by using
the fact that any point-line graph can be made naturally bipartite
by replacing every edge between a pair of points or a pair of lines
with a copy of the naturally bipartite graph $K_{3,3}$.

We close Section 4 by deriving a formula for the rank function of
$\scrm_{PL}$. When there are no points, this reduces to the formula
for the rank function of the cycle matroid of a graph. When there
are no lines it reduces to the formula for the rank function of the
2-dimensional bar-joint rigidity matroid given by Lov\'asz and
Yemini in \cite{LY}.

In Section 5 we consider the algorithmic implications of our
characterisation of $\scrm_{PL}$. We give a brief description of
Edmonds' algorithm \cite{E} for constructing a maximum independent
set in the union of two matroids. We then adapt the graph orientation
approach used by Berg and Jord\'an \cite{BJ} to give algorithms for augmenting independent sets in
$\scrm(\nu_L)$ and $\scrm(2\nu_P+\nu_L-2)$. Finally we describe
how these algorithms can be combined to determine the rank of
any generic set of point-line constraints.

\section{Definitions and preliminary results}

A {\em point-line graph} is a graph $G=(V,E)$ without loops together
with an ordered pair $(V_P,V_L)$ of, possibly empty, disjoint sets
whose union is $V$. We refer to vertices in $V_P$ and $V_L$ as {\em
point-vertices} and {\em line-vertices}, respectively. We label the
vertices as $V_P=\{u_1,\dots,u_{s}\}$ and $V_L=\{v_1,\dots,v_{t}\}$,
and the edges as $E=\{e_1,e_2,\ldots,e_m\}$. We use
$E_{PP}$, $E_{PL}$ and $E_{LL}$ to denote the  sets of edges
incident to two point-verices, to a point-vertex and a line-vertex and to
two line-vertices respectively. For $e \in E$, we write $e=xy$ to
mean that the end-vertices of $e$ are $x$ and $y$.
We will assume that graphs are simple
(without parallel edges) unless they are explicitly described as 
multigraphs. We supplement the above notation when it is
not obvious which graph we are referring to by using $V(G)$, $E(G)$,
etc.
We say that the point-line graph $G$ is {\em naturally bipartite} if
$G$ is a bipartite graph with bipartition $\{V_P,V_L\}$. In this
case $E_{PP}=E_{LL}=\emptyset$ and $E=E_{PL}$.

A {\em point-line framework} is a pair $(G,p)$ where $G$ is a
point-line graph
and
$p:V \to \real^{2}$. We put $p(u_i)=(x_i,y_i)$ for each $u_i\in V_P$
and $p(v_i)=(a_i,b_i)$ for each $v_i\in V_L$. This gives rise to a
geometric representation of $(G,p)$ by taking the point
corresponding to $u_i$ to have  cartesian coordinates $(x_i,y_i)$,
and the equation of the line corresponding to $v_i$ to be
$x=a_iy+b_i$.
We say that $(G,p)$ is {\em degenerate} if either $V_L=\emptyset$
and $p(u_i)=p(u_j)$ for all $u_i,u_j\in V_P$, or $V_P=\emptyset$ and
$a_i=a_j$ for all $v_i,v_j\in V_L$, and otherwise that $(G,p)$ is {\em non-degenerate}.




Given a point-line graph $G=(V,E)$, the {\em rigidity map} $f_G:\R^{2|V|}
\to \R^{|E|}$ is defined as follows. For
each $p\in \R^{2|V|}$ we consider the point-line framework $(G,p)$
and take $f_G(p)=(f_1(p),f_2(p),\ldots,f_m(p))$ where
$$
f_i(p)=\left\{
\begin{array}{ll}
(x_j-x_k)^2 + (y_j-y_k)^2 & \mbox{if $e_i=u_ju_k\in E_{PP}$}\\
(x_j-y_ja_k-b_k)(1+a_k^2)^{-\frac{1}{2}} & \mbox{if $e_i=u_jv_k\in E_{PL}$}\\
\tan^{-1}a_j-\tan^{-1}a_k & \mbox{if $e_i=v_jv_k\in E_{LL}$ and
$j<k$}
\end{array}
\right.
$$
The expressions for $f_i(p)$ are: the squared distance
between the points represented by $u_j,u_k$ when $e_i=u_ju_k\in
E_{PP}$; the signed distance between the point represented by $u_j$
and the line represented by $v_k$ when $e_i=u_jv_k\in E_{PL}$; the
angle between the lines represented by $v_j,v_k$ when $e_i=v_jv_k\in
E_{LL}$.

A point-line framework $(G,p)$ is {\em rigid} if there exists an
$\epsilon>0$ such that every point-line framework $(G,q)$ which
satisfies $f_G(q)=f_G(p)$ and $\|p(w)-q(w)\|<\epsilon$ for all $w\in
V$ can be obtained from $(G,p)$ by a rotation or translation of
$\R^2$.


Let $J(G,p)$ be the Jacobian matrix  of $f_G$  evaluated at some
point $p\in \R^{2|V|}$. Then $J(G,p)$ is a
$|E|\times 2|V|$-matrix with rows  indexed by $E$ and pairs of columns by $V$.
We label the two columns indexed by a vertex $u_i\in V_P$ as
$u_{i,x}$ and $u_{i,y}$, respectively, and the two columns indexed
by a vertex $v_i\in V_L$ as $v_{i,a}$ and $v_{i,b}$, respectively.
\begin{itemize}
\item A row in $J(G,p)$ indexed by an edge $e_i=u_ju_k \in E_{PP}$ has
entries
$$2(x_j-x_k),\  2(y_j-y_k),\ 2(x_k-x_j),\ 2(y_k-y_j) $$
in the columns indexed by $u_{j,x}$, $u_{j,y}$, $u_{k,x}$ and
$u_{k,y}$, respectively, and zeros elsewhere.
\item
A row in $J(G,p)$ indexed by an edge $e_i=u_jv_k \in E_{PL}$ has
entries
$$(1+a_k^2)^{-\frac{1}{2}},\ -a_k(1+a_k^2)^{-\frac{1}{2}},\ (-x_ja_k-y_j+a_kb_k)(1+a_k^2)^{-\frac{3}{2}},\ -(1+a_k^2)^{-\frac{1}{2}}$$
in the columns indexed by $u_{j,x}$, $u_{j,y}$, $v_{k,a}$ and
$v_{k,b}$, respectively, and zeros elsewhere.
\item
A row in $J(G,p)$ indexed by an edge $e_i=v_jv_k \in E_{LL}$ with
$j<k$ has entries
$$(1+a_j^2)^{-\frac{1}{2}},\ -(1+a_k^2)^{-\frac{1}{2}}$$
in the columns indexed by $v_{j,a}$ and $v_{k,a}$, respectively, and zeros
elsewhere.
\end{itemize}

The point $p$ is a {\em regular
point} of $f_G$ if the rank of $J(G,p)$
has its maximum value at $p$. Asimow and Roth \cite{AR} used methods from differential
geometry to show that the rigidity of a bar-joint framework $(G,p)$
is determined by the rank of $J(G,p)$ when $p$ is a regular point of
$f_G$. Similar arguments, based on the facts that the rotations and
translations of $\R^2$ generate a $3$-dimensional subspace of the
null space of $J(G,p)$ when $(G,p)$ is non-degenerate, and that
there exists an open neighbourhood $U$ of $p$ such that $\{q\in
U\,:\,f_G(q)=f_G(p)\}$
is a manifold of dimension $2|V|-\rank J(G,p)$
when $p$ is a regular point of $f_G$,
 can be
used to show:

\begin{lem}\label{lem:ar}
Let $(G,p)$ be a non-degenerate point-line framework.
Then\\
(a) $\rank J(G,p)\leq 2|V|-3$.\\
(b) If $\rank J(G,p)= 2|V|-3$ then $(G,p)$ is rigid.\\
(c) If $p$ is a regular point of $f_G$ and $\rank J(G,p)< 2|V|-3$
then $(G,p)$ is not rigid.
\end{lem}

It is straightforward to determine when a degenerate point-line
framework is rigid.  We will determine when a given point-line graph
$G$ can be realised as a non-degenerate point-line framework $(G,p)$
with $\rank J(G,p)= 2|V|-3$. To this end it will be helpful to apply
row and column operations to $J(G,p)$ to obtain the following
simpler $|E|\times 2|V|$-matrix $R(G,p)$.
\begin{itemize}
\item A row in $R(G,p)$ indexed by an edge $e_i=u_ju_k \in E_{PP}$ has
entries
$$x_j-x_k,\  y_j-y_k,\ x_k-x_j,\ y_k-y_j $$
in the columns indexed by $u_{j,x}$, $u_{j,y}$, $u_{k,x}$ and
$u_{k,y}$, respectively, and zeros elsewhere.
\item
A row in $R(G,p)$ indexed by an edge $e_i=u_jv_k \in E_{PL}$ has
entries
$$1,\ -a_k,\ -x_ja_k-y_j,\ -1$$
in the columns indexed by $u_{j,x}$, $u_{j,y}$, $v_{k,a}$ and
$v_{k,b}$, respectively, and zeros elsewhere.
\item
A row in $R(G,p)$ indexed by an edge $e_i=v_jv_k \in E_{LL}$ with
$j<k$ has entries
$$1,\ -1$$
in the columns indexed by $v_{j,a}$ and $v_{k,a}$, respectively, and zeros
elsewhere.
\end{itemize}

The Jacobian matrix $J(G,p)$ can be constructed from $R(G,p)$ using the
following operations. For each $v_i\in V_L$, we multiply the column
of $R(G,p)$ indexed by $v_{i,b}$ by $a_ib_i$ and
subtract it from the
column indexed by $v_{i,a}$,  then divide the resulting column by
$(1+a_i^2)$.
For each $e\in E_{PP}$, we multiply the row indexed by $e$ by $2$.
For each $e=u_jv_k\in E_{PL}$, we divide the row indexed by $e$ by
$(1+a_k^2)^{\frac{1}{2}}$. This construction immediately implies

\begin{lem}\label{lem:rigmatrix}
Let $(G,p)$ be a point-line framework. Then $\rank J(G,p) = \rank
R(G,p)$.
\end{lem}


We say that a point-line framework $(G,p)$ is: \textit{independent}
if $\rank R(G,p)=|E|$;  \textit{infinitesimally rigid} if it is
non-degenerate and $\rank R(G,p)=2|V|-3$; \textit{isostatic} if it
is both independent and infinitesimally rigid; {\em generic} if the
set of coordinates $\{x_i,y_i,a_j\,:\,u_i\in V_P\, ,v_j\in V_L\}$
are algebraically independent over $\rat$.

Note that the infinitesimal rigidity of $(G,p)$ is equivalent to the
condition that $\rank J(G,p)=2|V|-3$ whenever $(G,p)$ is
non-degenerate. Since the entries of $R(G,p)$ are polynomial functions
of the coordinates $x_i,y_i,a_j$, the rank of $R(G,p)$ will be
maximised whenever $(G,p)$ is generic. We say that the point-line
graph $G$ is {\em rigid} if some, or equivalently every, generic
realisation $(G,p)$ of $G$ is infinitesimally rigid. Lemmas
\ref{lem:ar} and \ref{lem:rigmatrix} tell us that $G$ is rigid if
and only if  some, or equivalently every, generic realisation
$(G,p)$ of $G$ is rigid.

The {\em rigidity matroid} $\scrm(G,p)$ of a point-line framework
 $(G,p)$ is the row matroid of its rigidity matrix $R(G,p)$. Its ground set
is $E$ and a set $S\subseteq E$ is independent if the rows of
$R(G,p)$ indexed by $S$ are linearly independent. The matroid
$\scrm(G,p)$ will be the same for all generic $(G,p)$. In this case
we refer to $\scrm(G,p)$ as the {\em rigidity matroid of $G$} and
denote it by $\scrm_{PL}(G)$, or simply $\scrm_{PL}$ when it is
obvious which graph we are referring to. We denote the rank  of
$\scrm_{PL}(G)$ by $r_{PL}(G)$. Thus $G$ is independent if and only
if $r_{PL}(G)=|E|$ and is rigid if and only if $r_{PL}(G)=2|V|-3$.

When $(G,p)$ is a point-line framework with  no point-vertices,
$R(G,p)$  is  the $\{0,1,-1\}$ edge/vertex incidence matrix of an
orientation of $G$. This immediately implies

\begin{lem}\label{lem:line_mat}
Let $(G,p)$ be a point-line framework with $V_P=\emptyset$. Then
$\scrm_{PL}(G,p)$ is the cycle matroid of $G$. In particular $(G,p)$
is independent if and only if $G$ is a forest.
\end{lem}

We can use this result to obtain the following necessary conditions
for $(G,p)$ to be independent.

\begin{lem} \label{easy-nec-con}
Suppose $(G,p)$ is an independent point-line framework and $H$ is a
subgraph of $G$ with $|E(H)|>0$. Then $|E(H)| \leq 2|V(H)|-3$. In
addition, if $V_P(H) = \emptyset$, then $|E(H)| \leq |V(H)|-1$.
\end{lem}
\bproof Since the rank of the rigidity matrix is maximised for
generic realisations of $G$, we may assume that $(G,p)$ is generic.
 The hypothesis that  $(G,p)$ is independent
implies that $(H,p|_H)$ is independent. We can now use Lemmas
\ref{lem:ar}(a) and \ref{lem:rigmatrix} to deduce that
$|E(H)|=\rank R(H,p|_H)=\rank J(H,p|_H)\leq 2|V(H)|-3$. When $V_P(H)
= \emptyset$, Lemma \ref{lem:line_mat} gives $|E(H)|=\rank
R(H,p|_H)\leq |V(H)|-1$. \eproof



\section{A count matroid for point-line graphs}

Given a point-line graph $G=(V,E)$,  the necessary conditions for
independence in $\scrm_\pl(G)$ given by Lemma \ref{easy-nec-con} are
not sufficient - we shall see that the point-line graph in Figure
\ref{fig1} is a counterexample. Indeed the family $\cal F$ of sets
satisfying these conditions need not even define a matroid on $E$.
On the other hand, we will show that $\scrm_\pl(G)$ is equal to the
matroid $\scrm_\sharp(G)$  on $E$ whose family of independent sets
is the maximum subset of $\cal F$ which satisfies the matroid
axioms.

We construct $\scrm_{\sharp}(G)$
from two simpler count matroids using
the operations of matroid union and Dilworth truncation.
We then use Lemma \ref{easy-nec-con} to prove the partial result
that every independent set in $\scrm_{\pl}(G)$ is an independent set
in $\scrm_\sharp(G)$. (More precisely we will show that if $I$ is an
independent set in some matroid on $E$ whose independent sets
satisfy Lemma \ref{easy-nec-con}, then $I$ is independent in
$\scrm_\sharp(G)$.) The reverse implication will be proved in the
next section.

We first recall some results from matroid theory. We refer a reader
unfamiliar with submodular functions and matroids to \cite{F}.

Let $E$ be a set. A {\em subpartition} of $E$ is a (possibly empty) collection of pairwise disjoint nonempty subsets of $E$.
A function $f:2^E\to \real$ is {\em nondecreasing} if $f(A)\leq f(B)$ for all
$A\subseteq B\subseteq E$,    {\em submodular} if $f(A)+f(B)\geq f(A\cup B)+f(A\cap B)$ for all
$A,B\subseteq E$,
and {\em intersecting submodular} if $f(A)+f(B)\geq f(A\cup B)+f(A\cap B)$ for all
$A,B\subseteq E$ with $A\cap B\neq \emptyset$.
We will need the following result of Dunstan \cite{Dunst}, see \cite[Theorem 12.1.1]{F}.

\begin{thm}\label{trunc} Suppose $E$ is a set and  $f:2^E\to \Z$ is intersecting submodular.
Let
$g:2^E\to \Z$ be defined by
putting $
g(\emptyset)=0$ and, for $\emptyset\neq A\subseteq E$,
$$g(A)=\min\left\{\sum_{i=1}^sf(A_i)\right\}$$
where the minimum is taken over all partitions $\{A_1,\ldots,A_s\}$ of $A$.
Then $g$ is submodular.
\end{thm}

\medskip

The following result of Edmonds \cite{E}, see \cite[Theorem 13.4.2]{F}, tells us how an
intersecting submodular function can be used to define a matroid on $E$.

\begin{thm}\label{induced} Suppose $E$ is a set and  $f:2^E\to \Z$ is nondecreasing, intersecting submodular and nonnegative on $2^E\sm \{\emptyset\}$.
Let
$$\scri=\{I\subseteq E\,:\,|J|\leq f(J) \mbox{ for all } \emptyset\neq J\subseteq I\}.$$
Then $\scri$ is the family of independent sets of a matroid $\scrm(f)=(E,\scri)$.
The rank of any $A\subseteq E$ in $\scrm(f)$ is given by
$$r(A)=\min\left\{\left|A\sm \bigcup_{i=1}^sA_i\right|+\sum_{i=1}^sf(A_i)\right\}$$
where the minimum is taken over all subpartitions $\{A_1,\ldots,A_s\}$ of $A$.
In addition we have:\\
(a) if $f(e)\leq 1$ for all $e\in E$, then
$r(A)=\min\left\{\sum_{i=1}^sf(A_i)\right\}$
where the minimum is taken over all partitions $\{A_1,\ldots,A_s\}$ of $A$;\\
(b) if $f$ is submodular and $f(\emptyset)=0$, then $r(A)=\min_{B\subseteq A}\left\{|A\sm B|+f(B)\right\}$;\\
(c) if $f$ is submodular, $f(\emptyset)=0$ and $f(e)\leq 1$ for all
$e\in E$, then $r(A)=f(A)$.
\end{thm}
The matroid $\scrm(f)$ given in Theorem \ref{induced} is referred to
as the matroid {\em induced by $f$}. The function $f-1$ will also
satisfy the hypotheses of Theorem \ref{induced} whenever $f$ has
$f(e)\geq 1$ for all $e\in E$. In this case we will refer to
$\scrm(f-1)$ as the {\em Dilworth truncation} of $\scrm(f)$.

The next result  describes a method to combine two matroids on the same ground set to obtain a new matroid.

\begin{thm}\cite{EF} \label{matunion}
Suppose that $\scrm_1=(E,\scri_1)$ and $\scrm_2=(E,\scri_2)$ are two matroids with the same ground set $E$.
Let
$$\scri=\{I_1\cup I_2\,:\,I_1\in \scri_1, I_2\in \scri_2\}.$$
Then $\scri$ is the family of independent sets of a matroid
$\scrm_1\vee\scrm_2$. The rank of any $A\subseteq E$ in
$\scrm_1\vee\scrm_2$ is given by
$$r(A)=\min_{B\subseteq A}\left\{r_1(B)+r_2(B)+\left|A\sm B\right|\right\}$$
where $r_1,r_2$ are the rank functions of $\scrm_1,\scrm_2$, respectively.
\end{thm}
The matroid $\scrm_1\vee\scrm_2$ is the {\em matroid union} of $\scrm_1$ and $\scrm_2$.

The expressions for the rank functions in Theorems \ref{induced}(b) and \ref{matunion} immediately give
 the following relationship between matroids induced by submodular functions and their matroid union.

\begin{lemma}\cite{PP}\label{submodunion}
Suppose $E$ is a set and  $f,g:2^E\to \Z$ are nondecreasing,
nonnegative submodular functions.
Then $f+g$ is a nondecreasing,
nonnegative submodular function and
$\scrm(f+g)=\scrm(f)\vee\scrm(g)$.
\end{lemma}
Note that the conclusion of Lemma \ref{submodunion} may not hold if the functions $f,g$ are allowed to take negative values on the empty set.





Let $G=(V,E)$ be a point-line graph. Recall that, for each $A\subseteq E$,
$\nu_P(A)$ and
$\nu_L(A)$ denote the numbers of point-vertices and line-vertices, respectively which are incident to edges in $A$.
It is easy to see that $\nu_P(A)$ and $\nu_L(A)$ are both nondecreasing, nonnegative, submodular functions on $2^E$.

\begin{lem}\label{distancerank} Let $G=(V,E)$ be a point-line graph.
Let $\rho:2^E\to \Z$ be defined by putting
$$\rho(A)=\min\left\{\sum_{i=1}^s (2\nu_P(A_i)+\nu_L(A_i)-2)\right\}$$
for all $A\subseteq E$, where
the minimum is taken over all partitions $\{A_1,\ldots,A_s\}$ of $A$.
Then $\rho$ and $\rho+\nu_L$ are  nondecreasing,  submodular and nonnegative, and $\rho+\nu_L-1$ is nondecreasing, submodular and  nonnegative on
$2^E\sm\{\emptyset\}$.
\end{lem}

\bproof
The function $2\nu_P+\nu_L-2$ is nondecreasing and submodular because $\nu_P$ and $\nu_L$ are both nondecreasing and submodular.
This, and Theorem \ref{trunc}, imply that $\rho$ is nondecreasing and submodular. The facts that $2\nu_P+\nu_L-2$ is nonnegative on
$2^E\sm\{\emptyset\}$ and $\rho(\emptyset)=0$ imply that $\rho$ is nonnegative. The assertion that  $\rho+\nu_L$ is  nondecreasing,  submodular and nonnegative now follows since $\nu_L$ is  nondecreasing,  submodular and nonnegative.
This immediately implies that $\rho+\nu_L-1$ is nondecreasing and submodular.
The assertion that $\rho+\nu_L-1$ is nonnegative on
$2^E\sm\{\emptyset\}$ follows since it is nondecreasing and $\rho(e)+\nu_L(e)-1=2\nu_P(e)+2\nu_L(e)-3= 1$ for all $e\in E$.
\eproof

We can now define our promised matroid $\scrm_\sharp(G)$ by putting
$\scrm_\sharp(G)=\scrm(\rho+\nu_L-1)$. Since $\scrm_\sharp(G)$ is
the Dilworth truncation of $\scrm(\rho+\nu_L)$, it will aid our
understanding of  $\scrm_\sharp(G)$ to express $\scrm(\rho+\nu_L)$
as a matroid union of two simpler matroids.

\begin{lem}\label{lem:rhomu_Lunion} Suppose $G=(V,E)$ is a point-line graph.
Then $$\scrm(\rho+\nu_L)=\scrm(2\nu_P+\nu_L-2)\vee\scrm(\nu_L).$$
\end{lem}
\bproof Since $\rho$ and $\nu_L$ are both nondecreasing and
submodular with $\rho(\emptyset)=0=\nu_L(\emptyset)$, we have
$\scrm(\rho+\nu_L)=\scrm(\rho)\vee\scrm(\nu_L)$ by Lemma
\ref{submodunion}. It remains to show that
$\scrm(\rho)=\scrm(2\nu_P+\nu_L-2)$. Let $r_1,r_2$ be the rank
functions of $\scrm(\rho)$ and $\scrm(2\nu_P+\nu_L-2)$,
respectively. We can use Theorem \ref{induced}(b) and the definition
of $\rho$ to deduce that
\begin{eqnarray*}
r_1(A)&=&\min_{B\subseteq A}\left\{\left|A\sm
B\right|+\rho(B)\right\}\\
&=& \min_{B\subseteq A}\min_{\scrp_B}\left\{\left|A\sm
B\right|+\sum_{B_i\in \scrp_B}(2\nu_P(B_i)+\nu_L(B_i)-2)\right\}
\end{eqnarray*}
for all $A\subseteq E$, where the second minimum runs over all
partitions $\scrp_B$ of $B$. On the other hand Theorem \ref{induced}
gives
$$r_2(A)=\min_{\scrq_A}\left\{\left|A\sm \bigcup_{B_i\in \scrq_A}B_i\right|+\sum_{B_i\in \scrq_A}
(2\nu_P(B_i)+\nu_L(B_i)-2)\right\}
$$
where the minimum runs over all subpartitions $\scrq_A$ of $A$. We
can now deduce that $r_1(A)=r_2(A)$ by putting $B=\bigcup_{B_i\in
\scrq_A} B_i$. \eproof

We will use Lemma \ref{lem:rhomu_Lunion} in Section \ref{sec:alg} to give an algorithm for constructing a maximal independent set in $\scrm(\rho+\nu_L)$. Our next two results will allow us to use this to find a maximal independent set in $\scrm(\rho+\nu_L-1)$.

Given a graph $G=(V,E)$ and $w_1,w_2\in V$ we denote the (multi)graph obtained by adding a new edge $w_1w_2$ to $G$ by $G+w_1w_2$. Note that if $e=w_1w_2$ already exists as an edge in $G$ then we add a new edge $e'$ parallel to $e$ in $G+w_1w_2$.

\begin{lem}\label{lem:covind} Suppose $G=(V,E)$ is a point-line graph.
Then the following statements are equivalent.\\
(a) $G$ is $\scrm(\rho+\nu_L-1)$-independent.\\
(b) $G+w_1w_2$ is $\scrm(\rho+\nu_L)$-independent for all $w_1,w_2\in V$.\\
(c) $G+w_1w_2$ is $\scrm(\rho+\nu_L)$-independent for all $w_1w_2\in E$.
\end{lem}
\bproof (a)$\Rightarrow$(b). Suppose $G$ is $\scrm(\rho+\nu_L-1)$-independent and let $w_1,w_2\in V$. Choose $A\subseteq E$. Then
$|A|\leq \rho(A)+\nu_L(A)-1$ and $|A+w_1w_2|=|A|+1\leq \rho(A)+\nu_L(A)\leq \rho(A+w_1w_2)+\nu_L(A+w_1w_2)$. Hence  $G+w_1w_2$ is $\scrm(\rho+\nu_L)$-independent.
\\[2mm]
(b)$\Rightarrow$(c) is immediate.
\\[2mm]
(c)$\Rightarrow$(a). We prove the contrapositive. Suppose $G$ is $\scrm(\rho+\nu_L-1)$-dependent. Then there exists a $\emptyset\neq B\subseteq E$ such that $|B|> \rho(B)+\nu_L(B)-1$. Choose a partition $\{B_1,B_2,\ldots B_t\}$ of $B$ such that $\rho(B)=\sum_{i=1}^t(2\nu_P(B_i)+\nu_L(B_i)-2)$. Choose $e=w_1w_2\in B_1$ and put $B_1'=B_1+w_1w_2$ and $B_i'=B_i$ for all $2\leq i\leq t$. Then
$\{B'_1,B'_2,\ldots B'_t\}$ is a partition of
$B+w_1w_2$, so $$\rho(B+w_1w_2)\leq \sum_{i=1}^t(2\nu_P(B_i')+\nu_L(B_i')-2)=\sum_{i=1}^t(2\nu_P(B_i)+\nu_L(B_i)-2)=\rho(B).$$
We also have $\nu_L(B+w_1w_2)=\nu_L(B)$. Hence
$|B+w_1w_2|=|B|+1 > \rho(B)+\nu_L(B)\geq \rho(B+w_1w_2)+\nu_L(B+w_1w_2)$ so $G+w_1w_2$ is
$\scrm(\rho+\nu_L)$-dependent.
\eproof

\begin{cor}\label{lem:covind1} Suppose $G=(V,E)$ is a point-line graph, $S \subset E$ is $\scrm(\rho+\nu_L-1)$-independent and $e\in E\sm S$. Then $S+e$ is $\scrm(\rho+\nu_L-1)$-independent if and only if $S+e+e^\prime$ is $\scrm(\rho+\nu_L)$-independent, where $e^\prime$ is a copy of $e$.
\end{cor}
\bproof If $S+e$ is $\scrm(\rho+\nu_L-1)$-independent then $S+e+e^\prime$ is $\scrm(\rho+\nu_L)$-independent by the equivalence of (a) and (c) in Lemma \ref{lem:covind}. Hence
suppose $S+e$ is $\scrm(\rho+\nu_L-1)$-dependent. Then there exists a $\emptyset\neq B\subseteq S+e$ such that $|B|> \rho(B)+\nu_L(B)-1$. We have $e \in  B$ since $S$ is $\scrm(\rho+\nu_L-1)$-dependent.
We can now show, as in the (c)$\Rightarrow$(a) part of the proof of Lemma  \ref{lem:covind}, that $\rho(B+e')\leq \rho(B)$ and    $\nu_L(B+e')=\nu_L(B)$.
This gives
$|B+e^\prime|=|B|+1 > \rho(B)+\nu_L(B)\geq \rho(B+e^\prime)+\nu_L(B+e^\prime)$ so $S+e+e^\prime$ is
$\scrm(\rho+\nu_L)$-dependent.

\eproof

The remainder of this section is devoted to showing that every independent set in $\scrm_\pl(G)$ is independent in $\scrm_\sharp(G)$.
The converse statement will be proved in the next section. Let $r_\pl(G)$ denote the rank of $\scrm_\pl(G)$.

We will need the following result from matroid theory.
Given a matroid $M=(E,r)$ and $F\subseteq E$
the {\em closure} of $F$ is cl$(F)=\{e\in E\,:\,r(F+e)=r(F)\}$.

\begin{lem} \label{cexchange}
Let $I_1,I_2,I_3$ be independent sets in a matroid $M$ with $I_1\cap
I_3=\emptyset$ and $I_2\subseteq $cl$(I_3)$. If $I_1\cup I_3$ is
independent then $I_1\cup I_2$ is independent.
\end{lem}
\bproof We prove the contrapositive. Assume that $I_1\cup I_2$ is
dependent. Let $C$ be a circuit in $I_1 \cup I_2 \cup I_3$ chosen
such that $C\cap (I_1\sm I_2)\neq \emptyset$ and, subject to this
condition, such that $C \cap (I_2 \setminus I_3)$ is minimal. (The
circuit $C$ exists since $I_1 \cup I_2$ is dependent and $I_2$ is
independent.) Choose $f\in C\cap (I_1\sm I_2)$. Suppose there exists
$e\in C \cap (I_2 \setminus I_3)$. Since $I_2\subset $cl$(I_3)$, $e$
is also contained in a circuit $C' \subseteq  I_3+ e$. Then
$f\not\in C'$ since $f\in I_1\sm I_2$, $e\in I_2$, and $I_1\cap
I_3=\emptyset$. By the matroid strong circuit exchange axiom, there
exists a circuit $C''$ with $f\in C'' \subseteq (C\cup C') - e
\subset I_1 \cup I_2 \cup I_3$. Then $C''$ contradicts the
minimality of $C\cap (I_2\sm I_3)$. Hence $C\subseteq I_1\cup I_3$ and so $I_1\cup
I_3$ is dependent.
\eproof

\begin{lem}\label{lem:partrank}
Let $G=(V,E)$ be a point-line graph and $\F=\{E_1,E_2,\ldots,E_t\}$ be a partition of $E$. Then
\begin{equation}\label{eq:partrank}
r_\pl(G)\leq \sum_{i=1}^t(2\nu_P(E_i)+\nu_L(E_i)-2)+|V_L|-c_L({\F})
\end{equation}
where $c_L({\F})$ is the number of components in the graph with
vertex set $\F$, in which two vertices $E_i,E_j$ are adjacent if
$V_L(E_i)\cap V_L(E_j)\neq \emptyset$.
\end{lem}
\bproof Let $H_i$ be the subgraph of $G$ induced by $E_i$. We may assume
that each $H_i$ is a complete graph since adding an edge between two
vertices of $H_i$ to $G$ will not change the right hand side of
(\ref{eq:partrank}) and cannot decrease the left hand side of
(\ref{eq:partrank}). Let $I$ be a maximal set of edges in $E_{LL}$
which are independent in $\scrm_\pl(G)$, $B$ be a base of
$\scrm_\pl(G)$ which contains $I$, $B_i=B\cap E_i$ and $I_i=B_i\cap
I$. Then
\begin{equation}\label{eq:part1}
r_\pl(G)=|B|=\sum_{i=1}^t|B_i|=\sum_{i=1}^t|B_i\sm I_i|+|I|\,.
\end{equation}

Lemma \ref{easy-nec-con} implies that $|I|\leq
|V_L|-c_L(\F')$ where $\F'=\{E_i\in \F\,:\,\nu_L(E_i)\neq 0\}$.
Since $c_L(\F)-c_L(\F')=|\F\sm \F'|$  inequality
(\ref{eq:partrank}) will follow from (\ref{eq:part1}) if we can show
that $|B_i\sm I_i|\leq 2\nu_P(E_i)+\nu_L(E_i)-2-\alpha_i$ for all $1\leq
i\leq t$, where $\alpha_i=1$ if $\nu_L(E_i)= 0$ and otherwise $\alpha_i=0$. This last inequality follows immediately from Lemma
\ref{easy-nec-con} when $\nu_L(E_i)=0$.
Hence we may suppose that $\nu_L(E_i)\neq 0$ and $\alpha_i=0$.




Let $I_i^*$ be a maximal set of edges in $E_i\cap E_{LL}$ which are
independent in $\scrm_\pl(H_i)$. Since $H_i$ is complete, $I_i^*$
will induce a spanning tree on the line-vertices of $H_i$, by Lemma
\ref{lem:line_mat}, and hence $|I_i^*|=\nu_L(E_i)-1$.  The
maximality of $I$ implies that $I_i^* \subseteq $cl$(I)$. Lemma \ref{cexchange}
(with $I_1 = B_i \setminus I_i$, $I_2 = I_i^*$, $I_3=I$) implies that
$(B_i \setminus I_i) \cup I_i^*$ is independent in $\scrm_\pl(G)$.

Lemma \ref{easy-nec-con} now gives $|(B_i\sm
I_i)\cup I_i^*|\leq 2\nu_p(E_i)+2\nu_L(E_i)-3$. Since
$|I_i^*|=\nu_L(E_i)-1$ we have $|B_i\sm I_i|\leq
2\nu_p(E_i)+\nu_L(E_i)-2$, as required. \eproof

\noindent {\bf Example} Consider the point-line graph $G=(V,E)$ on
the right hand side of Figure \ref{fig1}.  Let $E_1,E_2,E_3$ be the
sets of edges in the three copies of $K_4-e$.  Then
$\scrf=\{E_1,E_2,E_3\}$ is a partition of $E$. We have $\nu_L(E)=3$,
$c_L(\scrf)=1$ and $2\nu_P(E_i)+\nu_L(E_i)-2=4$ for all $1\leq i\leq
3$. Lemma \ref{lem:partrank} now gives $r_\pl(G)\leq 3+3\times 4
-1=14<2|V|-3$ so $G$ is not rigid.

\medskip

Lemma \ref{lem:partrank} and the definition of $\rho$ immediately give:

\begin{cor}\label{lem:pl_upperbound}
Let $G=(V,E)$ be a point-line graph. Then
$$r_\pl(G)\leq \rho(E)+|V_L|-1.$$
\end{cor}

\begin{lem}\label{lem:pl_ind}
Let $G=(V,E)$ be a point-line graph and $F\subseteq E$ be independent in $\scrm_\pl(G)$. Then
$F$ is independent in $\scrm_\sharp(G)$.
\end{lem}
\bproof Since $\scrm_\sharp(G)=\scrm(\rho+\nu_L-1)$, it will suffice
to show that $|F'|\leq \rho(F')+\nu_L(F')-1$ for all $F'\subseteq
F$. Let $H=G[F']$. Since $F'\subseteq F$, $F'$ is independent in
$\scrm_\pl(G)$ and hence $r_\pl(H)=|F'|$. We can now apply  Corollary
\ref{lem:pl_upperbound} to $H$ to deduce that $|F'|=r_\pl(H) \leq
\rho(F')+\nu_L(F')-1$. \eproof

\section{A characterisation of the generic point-line rigidity matroid}

We now complete the proof that $\scrm_\pl(G)=\scrm_\sharp(G)$.
Our approach is to first construct a linear representation for
$\scrm_\sharp(G)$ when $G$ is {\em naturally bipartite} i.e. $G$ is
a bipartite graph with bipartition $(V_P,V_L)$. We use this to
construct a point-line framework $(G,p)$ such that the rows of the
rigidity matrix $R(G,p)$ are linearly independent whenever $G$ is
$\scrm_\sharp$-independent and naturally bipartite. Together with
Lemma \ref{lem:pl_ind}, this will imply
$\scrm_\pl(G)=\scrm_\sharp(G)$ when $G$ is naturally bipartite. We
then deduce that this equality holds for an arbitrary point-line
graph by reducing to the bipartite case.

\subsection{Linear representations of point-line frames}

We will need the following result on linear representations of matroid unions. We include a proof for the sake of completeness.

\begin{lem}\cite[Lemma 7.6.14(1)]{B}\label{matunionrep}
Suppose that $\scrm_1$ and $\scrm_2$ are two matroids with the same ground set and that $\scrm_i$ is the row matroid of an $m\times n_i$ real matrix $M_i$ for $i=1,2$. Let $X$ be the $m\times m$ diagonal matrix $\diag(x_1,x_2,\ldots,x_m)$ where $x_1,x_2,\ldots,x_m$ are algebraically independent over $\rat[M_1]$.
Then $\scrm_1\vee\scrm_2$ is the row matroid of $M=(M_1,XM_2)$.
\end{lem}
\bproof Choose $F\subseteq E$. We show that $F$ is independent in $\scrm_1\vee\scrm_2$ if and only if the rows of $M$ labelled by $F$ are linearly independent.

We first suppose that $F$ is independent in $\scrm_1\vee\scrm_2$. Then there exists $I_i\in \scrm_i$ such that
$I_1\cap I_2=\emptyset$ and $I_1\cup I_2=F$. Let $x_i'=0$ for $i\in I_1$, $x_i'=1$ for $i\in E\sm I_1$, and put $M'=(M_1,X'M_2)$, where $X'=\diag(x'_1,x'_2,\ldots,x'_m)$. Reordering the rows of $M'$ if necessary, we have
$$M'[F]=\left(\begin{array}{cc}
M_1[I_1]&0\\
M_1[I_2]&M_2[I_2]
\end{array}\right).
$$
Since the rows of $M_1[I_1],M_2[I_2]$ are linearly independent, the rows of $M'[F]$ are linearly independent.
Since $x_1,x_2,\ldots,x_m$ are algebraically independent over $\rat[M_1]$, the rows of $M[F]$ are also linearly independent.

We next suppose that $F$ is dependent in $\scrm_1\vee\scrm_2$. Then there exists $K\subseteq F$ such that
$|F|>r(F)=r_1(K)+r_2(K)+|F\sm K|$ by Theorem \ref{matunion}. Reordering the rows of $M$ if necessary we have
$$M[F]=\left(\begin{array}{cc}
M_1[K]&XM_2[K]\\
M_1[F\sm K]&XM_2[F\sm K]
\end{array}\right).
$$
Since $\rank XM_2[K]=\rank M_2[K]$, this implies that
$$\rank M[F]\leq\rank M_1[K]+\rank M_2[K]+|F\sm K|=r_1(K)+r_2(K)+|F\sm K|<|F|$$
and hence the rows of $M[F]$ are linearly dependent.
\eproof

\medskip

A {\em point-line frame} is a triple $(G,t,c)$ where $G=(V,E)$ is a bipartite multigraph with bipartition
$V=V_P\cup V_L$, $t:V_L\to \R$ and $c:E\to \R$. We label the vertices in $V_P$ and $V_L$ as $u_1,u_2,...,u_m$
and $v_1,v_2,...,v_n$, respectively. We denote $t(v_j)$ by $t_j$ and $c(e)$ by $c_e$.
We associate the following three matrices with $(G,t,c)$.
\begin{itemize}
\item The {\em $A$-matrix}, $A(G,t,c)$, is the $|E|\times 2|V|$ matrix in which the entries in the row indexed by an edge $e=u_iv_j\in E$
are: $(1,t_j)$ in  the columns indexed by $u_i$; $(c_e,-1)$ in  the columns indexed by $v_j$; and zeros elsewhere.
\item The {\em $B$-matrix}, $B(G,c)$, is the $|E|\times (2|V_P|+|V_L|)$ matrix in which the entries in the row indexed by an edge $e=u_iv_j\in E$
are: $(1,c_e)$ in  the columns indexed by $u_i$; $-1$ in  the column indexed by $v_j$; zeros elsewhere.
\item The {\em $C$-matrix}, $C(G,t)$, is the $|E|\times (2|V_P|+|V_L|)$ matrix in which the entries in the row indexed by an edge $e=u_iv_j\in E$
are: $(1,t_j)$ in  the columns indexed by $u_i$; $-1$ in  the column indexed by $v_j$; zeros elsewhere.
\end{itemize}

The point-line frame $(G,t,c)$ is {\em generic} if the set $\{t_j,c_e\,:\,v_j\in V_L,e\in E\}$ is algebraically independent over $\rat$. We shall show that
 the matrices $B(G,c)$, $C(G,t)$ and $A(G,t,c)$  provide
linear representations  for the matroids $\scrm(2\nu_P+\nu_L-1)$,
$\scrm(2\nu_P+\nu_L-2)$ and $\scrm(2\nu_P+\nu_L-2)\vee
\scrm(\nu_L)$, respectively, when $(G,t,c)$ is generic. We first
need to  express $\scrm(2\nu_P+\nu_L-1)$ as the matroid union
$\scrm(\nu-1)\vee\scrm(\nu_P)$, where $\nu=\nu_P+\nu_L$ (and hence
$\scrm(\nu-1)$ is the well known {\em cycle matroid} of $G$). Note
that this does not follow from Lemma \ref{submodunion} because
$\nu(\emptyset)-1=-1$.

\begin{lem}\label{bipunion} Suppose $G=(V,E)$ is a naturally bipartite multigraph. Then
$$\scrm_G(2\nu_P+\nu_L-1)=\scrm_G(\nu-1)\vee\scrm_G(\nu_P).$$
\end{lem}
\bproof Let $r_1$ and $r_2$ be the rank functions of the matroids
$\scrm_G(\nu-1)$ and $\scrm_G(\nu_P)$, respectively. For $K\subset
E$, let $\Pi(K)$ be the partition of $K$ induced by the connected
components of $G[K]$. It is well known that $r_1(K)=\sum_{K_i\in
\Pi(K)}(\nu(K_i)-1)$ for all $F\subseteq E$,  and it is not
difficult to check that we also have
$r_2(K)=\sum_{K_i\in\Pi(K)}\nu_P(K_i)$. We can now apply \cite[Lemma
2.2]{KT} to deduce that
$\scrm_G(\nu-1)\vee\scrm_G(\nu_P)=\scrm_G(\nu+\nu_P-1)=\scrm_G(2\nu_P+\nu_L-1)$.
\eproof

\begin{lem}\label{lem:ABR-matrices} Suppose $(G,t,c)$ is a generic point-line frame. Then:\\
 (a) the row matroid of $B(G,c)$ is $\scrm_G(2\nu_P+\nu_L-1)$;\\
 (b) the row matroid of $C(G,t)$ is $\scrm_G(2\nu_P+\nu_L-2)$;\\
 (c) the row matroid of $A(G,t,c)$ is $\scrm_G(2\nu_P+\nu_L-2)\vee\scrm_G(\nu_L)$.
\end{lem}
\bproof (a) Let $\vec{G}$ be the directed graph obtained by
directing all edges of $G$ from $V_L$ to $V_P$. Then the matroid
$\scrm_G(\nu-1)$ is the row matroid of  the $|E|\times |V|$ matrix
$M_1$ which is the $\{0,1,-1\}$, edge/vertex incidence matrix for
$\vec{G}$. Similarly, the matroid $\scrm_G(\nu_P)$ is the row
matroid of  the $|E|\times |V_P|$ matrix $M_2$, which is the
$\{0,1\}$, edge/point-vertex incidence matrix for $G$. Lemma
\ref{matunionrep} now implies that
$\scrm_G(\nu-1)\vee\scrm_G(\nu_P)$ is the row matroid of $B(G,c)$,
after a suitable reordering of its columns, and (a) follows from
Lemma \ref{bipunion}.\\[2mm]
(b)
It will suffice to show that the rows of $C(G,t)$ are independent if and only if $E$ is independent in $\scrm_G(2\nu_P+\nu_L-2)$.

Suppose that the rows of $C(G,t)$ are independent. Choose
$F\subseteq E$ and let $H$ be the subgraph of $G$ induced by the
edges in $F$. Then the rows of $C(H,t)$ are independent, since the
entries in the rows of  $C(G,t)$ indexed by $F$ and columns indexed
by vertices not incident to $F$ are all zero. We can express the
vectors in $\Null C(H,t)$ in the form $(q,b)$ where $q:V_P(H)\to
\R^2$ and $b:V_L(H)\to \R$. Then the two vectors $(q,b)$ and
$(q',b')$  defined by $q_i=(1,0)$ and $q_i'=(0,1)$ for all $u_i\in
V_P(H)$, and  $b_j=1$ and $b_j'=t_j$ for all $v_j\in V_L(H)$, are
linearly independent and belong to $\Null C(H,t)$. This implies that
$$|F|=\rank C(H,t)\leq |V_L(H)|+2|V_P(H)|-2=\nu_L(F)+2\nu_P(F)-2.$$
Since this holds for all $F\subseteq E$, $E$ is independent in
$\scrm_G(2\nu_P+\nu_L-2)$.

We next suppose that $E$ is independent in
$\scrm_G(2\nu_P+\nu_L-2)$.
We may assume further that
each vertex in $V_P$ is incident with an edge of $G$ since, if this
is not the case, then we may add an edge incident to an isolated
vertex of $V_P$ in $G$ and preserve independence in
$\scrm_G(2\nu_P+\nu_L-2)$. Choose $e=u_iv_j\in E$ and let
$G+e'=(V,E+e')$ be the bipartite multigraph obtained by adding a new
edge $e'$ parallel to $e$ in $G$. Let $(G+e',\tilde t,\tilde c)$ be
a generic frame for $G+e'$ and $B(G+e',\tilde c)$ be its $B$-matrix.
The fact that $E$ is independent in $\scrm_G(2\nu_P+\nu_L-2)$
implies that $E+e'$ is independent in $\scrm_G(2\nu_P+\nu_L-1)$ and
hence, by (a), the rows of $B(G+e',\tilde c)$ are linearly
independent. We can represent each vector in the null space of
$B(G+e',\tilde c)$ as $(q,b)$ where $q:V_P\to \R^2$ and $b:V_L\to
\R$. Let $q(u_k)=(q_{k,1},q_{k,2})$ for all $u_k\in V_P$, and
$b(v_k)=b_k$ for all $v_k\in V_L$. Since $e=u_iv_j$ and $\rank
B(G,\tilde c)=\rank B(G+e',\tilde c)-1$, we can find a $(q,b)\in
\Null B(G,\tilde c)$ such that $q_{i,2}\neq 0$.

Repeating the above argument for each $e\in E$ and taking a suitable
linear combination of the vectors we obtain, we can construct a
$(q,b)\in \Null B(G,\tilde c)$ such that $q_{i,2}\neq 0$ for all
$u_i\in V_P$. The fact that $(q,b)\in \Null B(G,\tilde c)$ gives
\begin{equation}\label{enull1}
-b_j+q_{i,1}+\tilde c_e\,q_{i,2}=0\quad\mbox{for all}\quad e=u_iv_j\in E.
\end{equation}
Equation (\ref{enull1}) enables us to transform the matrix $C(G,b)$
to the matrix $B(G,\tilde c)$ by
subtracting $q_{i,1}$  times column $u_{i,1}$ from column $u_{i,2}$
and then dividing column $u_{i,2}$ by $q_{i,2}$, for all $u_i\in
V_P$. This implies that $\rank C(G,b)=\rank B(G,\tilde c)=|E|$.
Since $t$ is generic we have
$\rank C(G,t)\geq \rank C(G,b)$. Hence the rows of $C(G,t)$ are linearly independent.\\[2mm]
(c) This follows from (b), Lemma \ref{matunionrep} and the fact that the matroid $\scrm_G(\nu_L)$
is the row matroid of  the  $\{0,1\}$, edge/line-vertex incidence matrix for $G$.
\eproof

We could also have deduced Lemma \ref{lem:ABR-matrices}(b) from a result of Whiteley \cite[Theorem 4.1]{Wscene} or from the general theory of Dilworth truncations given in \cite{B}.

\begin{lem}\label{lem:bipind} Let $G=(V,E)$ be a naturally bipartite point-line graph.
Suppose that $E$ is independent in $\scrm_\sharp(G)$. Then $E$ is independent in $\scrm_\pl(G)$.
\end{lem}
\bproof Recall that $\scrm_\sharp(G)=\scrm(\rho+\nu_L-1)$ where $\rho$ is as defined in Lemma \ref{distancerank}.
We may assume that
each vertex in $V_L$ is incident with an edge of $G$ since, if this
is not the case, then we may add an edge incident to an isolated
vertex of $V_L$ in $G$ and preserve independence in
$\scrm_G(\rho+\nu_L-1)$. Choose $e=u_iv_j\in E$ and let
$G+e'=(V,E+e')$ be the bipartite multigraph obtained by adding a new
edge $e'$ parallel to $e$ in $G$. Let $(G+e', t, c)$ be a generic
frame for $G+e'$.
Since $E$ is independent in $\scrm_G(\rho+\nu_L-1)$, Lemma
\ref{lem:covind} implies that $E+e'$ is independent in
$\scrm_G(\rho+\nu_L)$ and hence, by Lemmas \ref{lem:rhomu_Lunion}
and \ref{lem:ABR-matrices}(c), the rows of $A(G+e',t, c)$ are
linearly independent. We can represent each vector in the null space
of $A(G+e',t, c)$ as $(q,h)$ where $q:V_P\to \R^2$ and $h:V_L\to
\R^2$. Let $q(u_k)=(q_{k,1},q_{k,2})$ for all $u_k\in V_P$ and
$h(v_k)=(h_{k,1},h_{k,2})$ for all $v_k\in V_L$. Since $e=u_iv_j$
and $\rank A(G,t, c)=\rank A(G+e',t, c)-1$, we can find a $(q,h)\in
\Null  A(G,t, c)$ such that $h_{j,1}\neq 0$.

Repeating this argument for each $e\in E$ and taking a suitable
linear combination of the vectors we obtain, we can construct a
$(q,h)\in \Null  A(G,t, c)$ such that $h_{j,1}\neq 0$ for all
$v_j\in V_L$. The fact that $(q,h)\in \Null  A(G,t, c)$ gives
\begin{equation}\label{enull2}
q_{i,1}+t_j\,q_{i,2}+ c_e\,h_{j,1}-h_{j,2}=0\quad\mbox{for all}\quad e=u_iv_j\in E.
\end{equation}
Construct a point-line framework $(G,p)$ by putting
$p(u_i)=(-q_{i,2},q_{i,1})$ for all $u_i\in V_P$ and
$p(v_j)=(-t_j,0)$ for all $v_j\in V_L$. Equation (\ref{enull2})
enables us to transform the rigidity matrix $R(G,p)$ to the
$A$-matrix $A(G,t,c)$ by
subtracting $h_{j,2}$  times column $v_{j,2}$ from column $v_{j,1}$
for all $v_j\in V_L$, and then dividing column $v_{j,1}$ by
$h_{j,1}$. This implies that $\rank R(G,p)=\rank A(G,t, c)=|E|$. It
follows that the rows of the rigidity matrix of any generic
realisation of $G$ as a point-line framework will be linearly
independent. Hence $E$ is independent in $\scrm_\pl(G)$. \eproof

\subsection{The non-bipartite case}

We reduce the general case to the naturally bipartite case by
replacing each `non-bipartite edge' by a copy of a naturally
bipartite $K_{3,3}$. We need the following lemma to show that this
operation preserves independence in $\scrm_\sharp$.

\begin{lem}\label{lem:cov_2sum} Let $G_1=(V_1,E_1)$ and $G_2=(V_2,E_2)$ be point-line graphs.
Suppose that $G_1$ and $G_2$  are $\scrm_\sharp$-independent,
$e=z_1z_2\in E_1$, and $V(G_1)\cap V(G_2)=\{z_1,z_2\}\in
\{V_P(G_1)\cap V_P(G_2),V_L(G_1)\cap V_L(G_2)\}$. Then
$G=(G_1-e)\cup G_2$ is  $\scrm_\sharp$-independent.
\end{lem}
\bproof Suppose $G=(V,E)$ is not $\scrm(\rho+\nu_L-1)$-independent.
Then there exists a nonempty $A\subseteq E$ such that
$\rho(A)+\nu_L(A)-1\leq |A|-1$. Since $G_1,G_2$ are
$\scrm(\rho+\nu_L-1)$-independent, we have $A\cap E_1\neq
\emptyset\neq A\cap E_2$. The definition of $\rho$ implies that
there exists a partition $\scrf=\{A_1,A_2,\ldots,A_t\}$ of $A$ such
that
\begin{equation}\label{eq:case1}
\sum_{i=1}^t (2\nu_P(A_i)+\nu_L(A_i)-2)+\nu_L(A)-1\leq |A|-1.
\end{equation}

We consider two cases.

\smallskip
\noindent
{Case 1:} $z_1,z_2\in V_L$.
We may assume that $\scrf$ has been chosen amongst all partitions satisfying (\ref{eq:case1}) to make $|\scrf|$ as large as possible. This choice ensures that $A_i\subseteq E_j$ for each $1\leq i\leq t$ and some $1\leq j\leq 2$, since if
$A_i\cap E_1\neq \emptyset\neq A_i\cap E_2$ then the partition $\scrf'=\scrf\sm \{A_i\})\cup\{A_i\cap E_1,A_i\cap E_2\}$ would still satisfy (\ref{eq:case1}).

For $1\leq j\leq 2$, let $\scrf_j=\{A_i\in \scrf\,:\,A_i\subseteq E_j\}$. Then $\scrf_j$ is a partition of $A_j'=A\cap E_j$. Since $V(A_1')\cap V(A_2')\subseteq \{z_1,z_2\}\subseteq V_L$, we have
\begin{eqnarray*}
|A_1'|+|A_2'|=|A|&\geq& \sum_{i=1}^t (2\nu_P(A_i)+\nu_L(A_i)-2)+\nu_L(A)\\
&\geq& \sum_{A_i\in \scrf_1} (2\nu_P(A_i)+\nu_L(A_i)-2)+\nu_L(A_1')-1+\\
&&\sum_{A_i\in \scrf_2} (2\nu_P(A_i)+\nu_L(A_i)-2)+\nu_L(A_2')-1.
\end{eqnarray*}
Since $G_j$ is $\scrm(\rho+\nu_L-1)$-independent we must have $|A_j'|=\sum_{A_i\in \scrf_j} (2\nu_P(A_i)+\nu_L(A_i)-2)+\nu_L(A_j')-1$ and $z_1,z_2\in V_L(A_j')$ for both $j=1, 2$. We can now put
$\scrf_1'=\scrf_1+\{e\}$. Then $\scrf_1'$ is a partition of $A_1''=A_1'+e$ which satisfies
\begin{align*}
\sum_{B_i\in \scrf_1'} &(2\nu_P(B_i)+\nu_L(B_i)-2)+\nu_L(A_1'')-1=\\
&\sum_{A_i\in \scrf_1} (2\nu_P(A_i)+\nu_L(A_i)-2)+\nu_L(A_1')-1
=|A_1'|=|A_1''|-1.
\end{align*}
This contradicts the hypothesis that $G_1$ is $\scrm(\rho+\nu_L-1)$-independent, since $A_1''\subseteq E_1$.

\smallskip
\noindent
{Case 2:} $z_1,z_2\in V_P$. We may assume that $\scrf$ has been
chosen amongst all partitions satisfying (\ref{eq:case1}) to make
$|\scrf|$ as small as possible. This choice ensures that $A_i\cap
A_j\cap V_P=\emptyset$ for all $1\leq i<j\leq t$, since if $A_i\cap
A_j\cap V_P\neq \emptyset$ then the partition $\scrf'=\scrf\sm
\{A_i,A_j\})\cup\{A_i\cup A_j \}$ would still satisfy
(\ref{eq:case1}). It follows, in particular that $z_1,z_2$ each
belong to at most one set $A_i\in \scrf$.

For $1\leq j\leq 2$, let $\scrf_j=\{A_i\cap E_j\,:\,A_i\in \scrf\mbox{ and } A_i\cap E_j\neq \emptyset\}$.
Then $\scrf_j$ is a partition of $A_j'=A\cap E_j$ for $j=1,2$. Since $V(A_1')\cap V(A_2')\subseteq \{z_1,z_2\}\subseteq V_P$, we have
\begin{eqnarray*}
|A_1'|+|A_2'|=|A|&\geq& \sum_{i=1}^t (2\nu_P(A_i)+\nu_L(A_i)-2)+\nu_L(A)\\
&\geq& \sum_{B_i\in \scrf_1} (2\nu_P(B_i)+\nu_L(B_i)-2)+\nu_L(A_1')+\\
&&\sum_{B_i\in \scrf_2} (2\nu_P(B_i)+\nu_L(B_i)-2)+\nu_L(A_2')-2s
\end{eqnarray*}
where $s=1$ if there exists $A_i\in \scrf$ with $z_1,z_2\in A_i$ and $A_i\cap E_1\neq \emptyset\neq A_i\cap E_2$, and $s=0$ otherwise.
Since $G_j$ is $\scrm(\rho+\nu_L-1)$-independent we must have $|A_j'|=\sum_{B_i\in \scrf_j} (2\nu_P(B_i)+\nu_L(B_i)-2)+\nu_L(A_j')-1$ for both $j=1,2$, and $s=1$. The construction of $\scrf_1$ now implies that $z_1,z_2\in V_P(B_k)$ for some $B_k\in \scrf_1$. We can now put
$\scrf_1'=\scrf_1\sm \{B_k\}\cup \{B_k+e\}$. Then $\scrf_1'$ is a partition of $A_1''=A_1'+e$ which satisfies
\begin{align*}
\sum_{C_i\in \scrf_1'} &(2\nu_P(C_i)+\nu_L(C_i)-2)+\nu_L(A_1'')-1=\\
&\sum_{B_i\in \scrf_1} (2\nu_P(B_i)+\nu_L(B_i)-2)+\nu_L(A_1')-1
=|A_1'|=|A_1''|-1.
\end{align*}
This contradicts the hypothesis that $G_1$ is $\scrm(\rho+\nu_L-1)$-independent, since $A_1''\subseteq E_1$.
\eproof

\begin{lem}\label{lem:K33} Let $G=(V,E)$ be a naturally bipartite point-line graph which  is
isomorphic to $K_{3,3}$. Then $G$ is $\scrm_{\sharp}$-independent.
\end{lem}
\bproof Choose $e=uv\in E$. Let $B$ be a perfect matching of $G+uv$ which contains $uv$ and put $A=E(G+uv)\sm B$. Then
$E(G+uv)=A\cup B$, $A$ is $\scrm(2\nu_P+\nu_L-2)$-independent
and $B$ is $\scrm(\nu_L)$-independent. Symmetry and Lemmas
\ref{lem:rhomu_Lunion} and \ref{lem:covind} now imply that $G$ is
$\scrm_{\sharp}$-independent.
\eproof

\begin{lem}\label{lem:ind} Let $G=(V,E)$ be a point-line graph which is $\scrm_{\sharp}$-independent.
Then $G$ is $\scrm_\pl$-independent.
\end{lem}
\bproof We use induction on $|E_{PP}\cup E_{LL}|$. The lemma follows
from Lemma \ref{lem:bipind} when  $E_{PP}\cup E_{LL}=\emptyset$.
Hence we may suppose that we have an edge $e=w_1w_2\in E_{PP}\cup
E_{LL}$. Let $H$ be a naturally bipartite point-line graph which is
isomorphic to $K_{3,3}$ and label its vertices such that
$\{w_1,w_2\}=V_P(G)\cap V_P(H)$ or $\{w_1,w_2\}=V_L(G)\cap V_L(H)$.
Let $G^+=(G-e)\cup H$. Then Lemmas \ref{lem:cov_2sum} and
\ref{lem:K33} imply that $G^+$ is $\scrm_{\sharp}$-independent. We
can now use induction to deduce that $G^+$ is
$\scrm_\pl$-independent. This implies that $G$ is
$\scrm_\pl$-independent (since if $G$ were $\scrm_\pl$-dependent
then the matroid circuit axiom and the fact that $H+e$ is
$\scrm_\pl$-dependent would imply that $G^+=(G\cup (H+e))-e$
contains an $\scrm_\pl$-circuit). \eproof

\begin{thm}\label{thm:equal} Let $G=(V,E)$ be a point-line graph.
Then $\scrm_\pl(G)=\scrm_\sharp(G)$.
\end{thm}
\bproof Choose $F\subseteq E$. If $F$ is independent in $\scrm_\pl(G)$ then $F$ is independent in $\scrm_\sharp(G)$ by Lemma \ref{lem:pl_ind}. On the other hand,
if $F$ is independent in $\scrm_\sharp(G)$ then we may apply Lemma \ref{lem:ind} to $G[F]$ to deduce that $F$ is independent in $\scrm_\pl(G)$.
\eproof

\subsection{The rank function}

Theorem \ref{thm:equal} tells us that $\scrm_\pl(G)$ is the matroid
induced by $\rho+\nu_L-1$ and we can now use Theorem
\ref{induced}(a) to deduce that its rank function is given by
\begin{equation}\label{eq:rank1}
r_\pl(A)=\min_\scrf \left\{\sum_{A_i\in \scrf}(\rho(A_i)+\nu_L(A_i)-1)\right\}
\end{equation}
for all $A\subseteq E$, where the minimum is taken over all partitions $\scrf$ of $A$.
We close this section by obtaining an  expression for $r_\pl(A)$ in terms of $\nu_P$, $\nu_L$ and the function $c_L$ defined in the statement of Lemma \ref{lem:partrank}. 

\begin{thm}\label{thm:rank} Let $G=(V,E)$ be a point-line graph and $A\subseteq E$. Then
\begin{equation}\label{eq:rank0}
 r_\pl(A)=\nu_L(A)+\min_\scrf \left\{\sum_{A_i\in\scrf}(2\nu_P(A_i)+\nu_L(A_i)-2)-c_L(\scrf)\right\}
 \end{equation}
 where the
minimum is taken over all partitions $\scrf$ of $A$.
\end{thm}
\bproof The fact that the right hand side gives an upper bound on
$r_\pl(A)$ follows from Lemma \ref{lem:partrank}.
Hence it will suffice to show that
there exists a partition of $A$ which gives equality in
(\ref{eq:rank0}).

By (\ref{eq:rank1}), we may choose a partition $\scrp$ of $A$ such
that \begin{equation}\label{eq:rank2} r_\pl(A)=\sum_{B_i\in
\scrp}(\rho(B_i)+\nu_L(B_i)-1)
 \end{equation}
 and, subject to this condition, $|\scrp|$ is as small as possible. We claim that $V_L(B_i)\cap V_L(B_j)= \emptyset$ for all distinct $B_i,B_j\in\scrp$.

Suppose to the contrary that $V_L(B_i)\cap V_L(B_j)\neq \emptyset$. Let $\scrq=\scrp\sm\{B_i,B_j\}\cup\{B_i\cup B_j\}$. Then
\begin{align*}
\sum_{B_k\in \scrp}&(\rho(B_k)+\nu_L(B_k)-1)-\sum_{C_k\in \scrq}(\rho(C_k)+\nu_L(C_k)-1)=\\
&\rho(B_i)+\rho(B_j)-\rho(B_i\cup B_j)+\nu_L(B_i)+\nu_L(B_j)-\nu_L(B_i\cup B_j)-1\geq 0
\end{align*}
since $\rho,\nu_L$ are submodular and nonnegative, and $\nu_L(B_i\cap B_j)\geq 1$. Equations (\ref{eq:rank1}) and (\ref{eq:rank2}) now imply that $r_\pl(A)=\sum_{C_k\in \scrq}(\rho(C_k)+\nu_L(C_k)-1)$ and hence $\scrq$ contradicts the choice of $\scrp$.

The fact that the sets in $\scrp$ are line-vertex disjoint now gives
\begin{equation}\label{eq:rank3}
r_\pl(A)=\sum_{B_i\in \scrp}(\rho(B_i)+\nu_L(B_i)-1)=
\nu_L(A)+\left(\sum_{B_i\in \scrp}\rho(B_i)\right)-|\scrp|\,.
\end{equation}

For each $B_i\in \scrp$ we have
$$\rho(B_i)=\min_{\scrp_i}\left\{\sum_{A_{i,j}\in \scrp_i}(2\nu_P(A_{i,j})+\nu_L(A_{i,j})-2)\right\},$$
where the minimum is taken over all partitions $\scrp_i$ of $B_i$.
Choose a partition $\scrf_i$ of $B_i$ such that
$\rho(B_i)=\sum_{A_{i,j}\in
\scrf_i}(2\nu_P(A_{i,j})+\nu_L(A_{i,j})-2)$ and, subject to this
condition, such that $|\scrf_i|$ is as small as possible.
This choice ensures that $A_{i,j}\cap
A_{i,k}\cap V_P=\emptyset$ for all distinct $A_{i,j},A_{i,k}\in \scrf_i$, since if $A_{i,j}\cap
A_{i,k}\cap V_P\neq \emptyset$ then the partition $\scrf_i'=\scrf\sm
\{A_{i,j},A_{i,k}\})\cup\{A_{i,j}\cup A_{i,k} \}$ would contradict the minimality of $|\scrf_i|$.

 Let $\scrf=\bigcup_{B_i\in \scrp}\scrf_i$.
Then $\scrf$ is a partition of $A$ and we may use (\ref{eq:rank3})
to obtain
$$r_\pl(A)=\nu_L(A)+\sum_{A_{i,j}\in \scrf}(2\nu_P(A_{i,j})+\nu_L(A_{i,j})-2)-|\scrp|.
 $$
 It remains to show that $|\scrp|=c_L(\scrf)$, or equivalently, that $c_L(\scrf_i)=1$ for all $B_i\in \scrp$.

Suppose to the contrary that $c_L(\scrf_i)\geq 2$ for some $B_i\in \scrp$. Then there exists a partition
of $\scrf_i$ into two sets $\scrf_i',\scrf_i''$ such that $V_L(A_{i,j})\cap V_L(A_{i,k})=\emptyset$ for all $A_{i,j}\in \scrf_i'$ and $A_{i,k}\in \scrf_i''$.
Let $B_i'=\bigcup_{A_{i,j}\in \scrf_i'}A_{i,j}$ and
$B_i''=\bigcup_{A_{i,k}\in \scrf_i''}A_{i,k}$.
Then
 \begin{align}
\rho(B_i)&=\sum_{A_{i,j}\in \scrf_i}(2\nu_P(A_{i,j})+\nu_L(A_{i,j})-2)\nonumber\\
&=\sum_{A_{i,j}\in \scrf_i'}(2\nu_P(A_{i,j})+\nu_L(A_{i,j})-2)+\sum_{A_{i,k}\in \scrf_i''}(2\nu_P(A_{i,k})+\nu_L(A_{i,k})-2)\nonumber\\
&\geq \rho(B_i')+\rho(B_i'').\label{eq:rank4}
\end{align}
Let $\scrr=(\scrp\sm\{B_i\})\cup\{B_i',B_i''\}$. Then $\scrr$ is a partition of $A$ into line-vertex disjoint sets and
\begin{align*}
\sum_{B_j\in \scrp}&(\rho(B_j)+\nu_L(B_j)-1)-\sum_{C_j\in \scrr}(\rho(C_j)+\nu_L(C_j)-1)\\
&=
\nu_L(A)+\left(\sum_{B_j\in \scrp}\rho(B_j)\right)-|\scrp|-\nu_L(A)-\left(\sum_{C_j\in \scrr}\rho(C_j)\right)+|\scrr|\\
&=
\rho(B_i)-\rho(B_i')-\rho(B_i'')+1\geq 1
\end{align*}
by (\ref{eq:rank4}). This contradicts the minimality of $\sum_{B_j\in \scrp}(\rho(B_j)+\nu_L(B_j)-1)$. Hence $c_L(\scrf_i)=1$.
\eproof

Note that the expression for $r_\pl(A)$ in Theorem \ref{thm:rank} reduces to the Lov\'asz-Yemini rank formula  for the bar-joint rigidity matroid \cite{LY} when $V_L=\emptyset$ since we have $\nu_L(A)=0=v_L(A_i)$ and $c_L(\scrf)=|\scrf|$.



\section{Algorithmic Implications}\label{sec:alg}
We will describe an efficient (i.e. polynomial) algorithm to determine the rank $r_{PL}(G)$ of a point-line graph $G=(V,E)$.
By basic facts from matroid theory, $r_{PL}(G)$ is the size of any maximal independent set $I \subseteq E$.
Hence, by Theorem \ref{thm:equal}, we will have an efficient algorithm to determine
$r_{PL}(G)$ provided we can efficiently determine if $I+e$
is $\scrm_\sharp$-independent when $I$ is $\scrm_\sharp$-independent and $e\in E\sm I$.
By definition $\scrm_\sharp = \scrm(\rho+\nu_L-1)$, and Corollary \ref{lem:covind1} implies that  $I+e$ is $\scrm(\rho+\nu_L-1)$-independent if and only if $I+e+e^\prime$ is $\scrm(\rho+\nu_L)$-independent, where $e^\prime$ is a copy of $e$. Hence we require an efficient algorithm to test if $I+e+e^\prime$ is $\scrm(\rho+\nu_L)$-independent. 

Lemma \ref{lem:rhomu_Lunion} tells us that
$\scrm(\rho+\nu_L)=\scrm(2\nu_P+\nu_L-2)\vee \scrm(\nu_L)$. This
implies that we may use Edmonds' algorithm for matroid union
\cite{Ematunion}, see also Gabow and Westermann \cite{GW}, to
efficiently test for independence in $\scrm(\rho+\nu_L)$, as long as
we can efficiently test for independence in both
$\scrm(2\nu_P+\nu_L-2)$ and $\scrm(\nu_L)$. The latter can be
accomplished by using existing algorithms for count matroids
based on network flows \cite{I,S}, graph orientations \cite{BJ,F} or the pebble game
\cite{JH,LS}.

\medskip

\subsection{Matroid union and Augmenting Paths}\label{subsec:union}
We now give a brief description of Edmonds' algorithm for matroid
union. We refer the reader to \cite{GW} for more details. We are
given matroids $\scrm_1, \scrm_2$ with the same groundset $E$, an
independent set $I$ of $\scrm_1\vee \scrm_2$, a partition
$(I_1,I_2)$ of $I$ with $I_q$ independent in $\scrm_q$ for $q\in
\{1,2\}$. For any $e\in E\sm I_q$ such that $I_q+e$ is not
independent in $\scrm_q$ let $C(I_q,e,\scrm_q)$ be the unique
circuit of $\scrm_q$ contained  in $I_q+e$. We determine whether
$I+e$ is independent in $\scrm_1\vee \scrm_2$ by searching for an
{\em augmenting path}. This is a sequence of elements
$e=e_0,e_1,\ldots,e_s$ of $I+e$ with the following properties (where
subscripts $j$ on $I_j$ and $\scrm_j$ are to be read modulo two).

For some $q\in \{1,2\}$
\begin{itemize}
\item $e_{i+1}\in I_{q+i}$ for all $0\leq i\leq s-1$.
\item For all $0\leq i\leq s-1$, $I_{q+i}+e_i$ is dependent in
$\scrm_{q+i}$ and $e_{i+1}\in C(I_{q+i},e_i,\scrm_{q+i})-e_i$.
\item For all $1\leq i+1<j\leq  s$, $e_{j}\not\in C(I_{q+i},e_i,\scrm_{q+i})$.
\item $I_{q+s}+e_s$ is independent in $\scrm_{q+s}$.
\end{itemize}

If we find an augmenting path then we conclude that $I+e$ is
independent in $\scrm_1\vee \scrm_2$. We output $I+e$ together with
the partition $(I_q',I_{q+1}')$ of $I+e$ where $I_q'=I_q\triangle
\{e_0,e_1,\ldots, e_s\}$ is independent in $\scrm_q$,
$I_{q+1}'=I_{q+1}\triangle \{e_1,\ldots, e_s\}$ is independent in in
$\scrm_{q+1}$, and $\triangle$ denotes symmetric difference. The
requirement that the augmenting path has no short cuts (third bullet
point above) ensures that $I_q'$ is independent in $\scrm_q$ for
$q\in \{1,2\}$. If no augmenting path exists then we conclude that
$I+e$ is dependent in $\scrm_1\vee \scrm_2$.


To implement this algorithm
we need subroutines which determine whether $I_q+e$ is independent
in $\scrm_q$, and determine the unique circuit of $\scrm_q$
contained in $I_q+e$ when it is dependent.  We will adapt the algorithms for count matroids  given by Berg and
Jord\'an \cite{BJ} to obtain these
subroutines for $\scrm_1=\scrm(2\nu_P+\nu_L-2)$ and
$\scrm_2=\scrm(\nu_L)$.

\subsection{Graph orientations}

Let $G=(V,E)$ be a graph. An {\em orientation  of $G$} is  a
directed graph $D$ obtained by replacing each edge $wz\in E$ by a
directed edge (directed from $w$ to $z$ or from $z$ to $w$).
For $w \in V$, let $d^-_D(w)$ be the number of edges entering $w$ in $D$. 
Given a map $g:V \to \mathbb Z^+=\{n\in \mathbb Z\,:\,n\geq 0\}$, 
an orientation $D$ of $G$ is said to be a {\em $g$-orientation}  if $d^-_D(w) \leq g(w)$ for all $w \in V$.
For $X \subseteq V$, let $g(X)=\sum_{v \in X}g(v)$
 and let $i_G(X)$ be the number of edges
 induced by $X$.
The following result
 of Frank and Gy\'arf\'as \cite{FG} characterises when a graph has a $g$-orientation.

\begin{thm} \label{FandG}
Let $G=(V,E)$ be a graph and $g:V\to \mathbb Z^+$. Then  $G$ has a $g$-orientation if and only if $i_G(X) \leq g(X)\ for\ all\ X \subseteq V.$
\end{thm}


Given a point-line graph $G=(V,E)$ and $i,j\in \mathbb Z^+$, let
$g_{i,j}:V \to \mathbb Z^+$ be defined by $g_{i,j}(x)=i$  if $x \in V_P$ and $g_{i,j}(x)=j$ if $x \in V_L$.
We will be interested in (special kinds of) $g_{2,1}$- and $g_{0,1}$-orientations but it will be more efficient to consider general $g_{i,j}$-orientations.
Given $w,z\in V$ and $k\in \Z^+$,
a {\em $g_{i,j,k}^{wz}$-orientation of $G$} is a $g_{i,j}$-orientation $D$ such that $$d^-_D(w)+d^-_D(z) \leq g_{i,j}(w)+g_{i,j}(z)-k.$$

We will assume henceforth that $k\leq \min\{2i,2j\}$. In this case,
the nondecreasing submodular function $i\nu_P+j\nu_L-k$ is
nonnegative on $2^E\sm \{\emptyset\}$ and hence induces a matroid on
$E$ by Theorem \ref{induced}. Let
$\scrm(i,j,k)=\scrm(i\nu_P+j\nu_L-k)$. The following results give
relationships between $\scrm(i,j,k)$ and $g_{i,j}$-orientations.

%
%
%
%

\begin{lem} \label{lem:inddir}
Suppose $G=(V,E)$ is a point-line graph. Then $G$ is $\scrm(i,j,k)$-independent if and only if $G$ has a $g_{i,j,k}^{wz}$-orientation for all $w,z \in V$.
\end{lem}
\bproof
Suppose that
$G$ has no $g_{i,j,k}^{wz}$-orientation for some $w,z \in V$. Let $h:V\to \mathbb Z^+$ be such that $h(w)+h(z)=g_{i,j}(w)+g_{i,j}(z)-k$,
$h(x)\leq g_{i,j}(x)$ for all $x \in\{w,z\}$, and $h(x)=g_{i,j}(x)$ for all $x \in V \setminus \{w,z\}$. Then $G$ has no $h$-orientation because any $h$-orientation would be a $g_{i,j,k}^{wz}$-orientation. By Lemma \ref{FandG} there is a set $X \subseteq V$ with $i_G(X) > h(X)$. Let $A \subseteq E$ be the set of edges of $G$ induced by $X$. We have $|A|=i_G(X) > h(X) \geq g_{i,j}(X)-k=i\nu_P(A)+i\nu_L(A)-k$.
Hence $G$ is not $\scrm(i,j,k)$-independent.

Conversely suppose that $G$ is not $\scrm(i,j,k)$-independent. Then  $|A| > i\nu_P(A)+j\nu_L(A)-k$ for some nonempty $A \subseteq E$. Let $X=V(A)$.
Then $i_G(X)\geq |A| > i\nu_P(A)+j\nu_L(A)-k$ and Lemma \ref{FandG} implies that $G$ has no $g_{i,j,k}^{wz}$-orientation for any distinct $w,z \in V(A)$.
\eproof

For $I\subseteq E$ let $G(I)$ denote the spanning subgraph of $G$ with edge set $I$.


\begin{lem}\label{lem:extendI}
Let $G=(V,E)$ be a point-line graph, $I \subset E$ be independent in $\scrm(i,j,k)$ and $e=wz \in E \setminus I$.
 Then $I+e$ is  independent in $\scrm(i,j,k)$ if and only if $G(I)$ has a $g_{i,j,k+1}^{wz}$-orientation.
\end{lem}
\bproof
Suppose that $I+e$ is not independent in $\scrm(i,j,k)$. Then there exists an $A \subseteq I+e$ such that $|A| > i\nu_P(A)+j\nu_L(A)-k$.  Since $I$ is independent in $\scrm(i,j,k)$, we must have  $e \in A$, and $|A-e| = i\nu_P(A)+j\nu_L(A)-k$. Let  $X=V(A)$. Then $w,z \in X$ and
$i_{G(A-e)}(X)\geq |A-e| = i\nu_P(A)+j\nu_L(A)-k$.
Lemma \ref{FandG}
now implies that $G(I)$ has
no $g_{i,j,k+1}^{wz}$-orientation.

Conversely, suppose that  $G(I)$ has no $g_{i,j,k+1}^{wz}$-orientation. Then $G(I+e)$ has no $g_{i,j,k}^{wz}$-orientation.
Lemma \ref{lem:inddir} now implies that $I+e$ is dependent in $\scrm(i,j,k)$.
\eproof

Suppose $\scrm$ is a matroid, $I$ is an independent set in $\scrm$
and $e$ is an element of $\scrm$ such that $I+e$ is dependent.
Constructing the circuit $C(I,e,\scrm)$ is an important step
in the algorithm for matroid union outlined in Section
\ref{subsec:union}.
Our next lemma
tells us how to do this when $\scrm=\scrm(i,j,k)$.

\begin{lem}\label{lem:fundcirc}
Let $G=(V,E)$ be a point-line graph, $I \subset E$ be independent in $\scrm(i,j,k)$ and $I+e$ be dependent for some $e=wz \in E \setminus I$.
 Then $G(I)$ has a $g_{i,j,k}^{wz}$-orientation $D$. Furthermore, if $Y$ is the set of all vertices which are connected to $\{w,z\}$ by directed paths in $D$, and $F$ is the set of all edges of $I$ which are induced by $Y$, then $C(I,e)=F+e$.
\end{lem}
\bproof
The fact that $G(I)$ has a $g_{i,j,k}^{wz}$-orientation follows from Lemma \ref{lem:inddir}.
By definition, $C=C(I,e)$ is the minimal $\scrm(i,j,k)$-dependent subset of $I+e$. Hence $|C|=i\nu_P(C)+j\nu_L(C)-k+1$ and
$|C-e|=i\nu_P(C-e)+j\nu_L(C-e)-k$. Let $Y'=V(C)$. Since $I$ is independent, $e\in C$ and hence $w,z\in Y'$. Let $G'=(Y',C-e)$ and $D'$ be the restriction of
$D$ to $G'$. The facts that $D'$ is a
$g_{i,j,k}^{wz}$-orientation of $G'$ and $|C-e|=i\nu_P(C-e)+j\nu_L(C-e)-k$ imply that $d_{D'}^-(y)=g_{i,j}(y)$ for all $y\in Y\sm \{w,z\}$ and
$d_{D'}^-(w)+d_{D'}^-(z)=g_{i,j}(w)+g_{i,j}(z)-k$. Since $D$ is a $g_{i,j,k}^{wz}$-orientation of $G$, this gives $d_{D'}^-(y)=d_{D}^-(y)$ for all
$y\in Y'$. Thus there are no directed edges in $D$ from $V\sm Y'$ to $Y'$ and hence $Y\subseteq Y'$.

On the other hand, the definition of $Y$ implies that there are no directed edges in $D$ from $V\sm Y$ to $Y$. Thus, if $F$ is the set of edges of $I$ induced by $Y$ and $D''$ is the restriction of $D$ to $(Y,F)$, we will have $d_{D''}^-(y)=d_{D'}^-(y)$ for all $y\in Y$. This gives
$$|F|=\sum_{y\in Y}d_{D''}^-(y)=\sum_{y\in Y}d_{D'}^-(y)=i\nu_P(F)+j\nu_L(F)-k.$$
This implies that $F+e$ is dependent in $\scrm(i,j,k)$ and the minimality of $C$ now gives $C=F+e$ (and $Y'=Y$).
\eproof

Lemmas \ref{lem:inddir}, \ref{lem:extendI} and \ref{lem:fundcirc}
give rise to the following linear time algorithm which either
increases the size of an independent set $I$ in $\scrm(i,j,k)$ by
adding a new element $e$ to it,  or finds the fundamental circuit
$C(I,e)$ when $I+e$ is dependent.

Suppose that we are given $I$ together with a $g_{i,j}$-orientation
$D$ of $G(I)$ and an edge $e=wz\in E\sm I$. If $D$ is a
$g_{i,j,k+1}^{wz}$-orientation of $G(I)$ then we conclude that $I+e$
is independent. We orient $e$ so that $D+e$ is a
$g_{i,j,k}^{wz}$-orientation of $G(I+e)$ and output $(I+e,D+e)$.
Otherwise we construct the set $Y$ of all vertices which are
connected to $\{w,z\}$ by directed paths in $D$. If some vertex
$y\in Y\sm \{w,z\}$ has $d_D^-(y)<g_{i,j}(y)$ then we construct a
new orientation by reversing the direction of all edges on a
directed path from $y$ to $\{w,z\}$ in $D$ and then iterate. (Note
that the reorientation will reduce $d_D^-(w)+d_D^-(z)$ by one.)
After at most $k$ iterations, we will arrive at either a
$g_{i,j,k+1}^{wz}$-orientation of $G(I)$, or a
$g_{i,j,k}^{wz}$-orientation with  $d_D^-(y)=g_{i,j}(y)$ for all
$y\in Y\sm \{w,z\}$. In the latter case we conclude that $I+e$ is
dependent and output  $C(I,e)=F+e$, where $F$ is the set of all
edges of $I$ induced by $Y$.

We may combine this algorithm for $\scrm(2,1,2)$ and $\scrm(0,1,0)$
with the augmenting path algorithm of Section \ref{subsec:union} to
give an O$(|V|^2)$-algorithm which determines whether $I+e$ is
independent in $\scrm_\pl$ and hence obtain an
O$(|V|^2|E|)$-algorithm for constructing a maximum independent set
in $\scrm_\pl$. A more detailed description of this algorithm as a
pebble game is given in the arxiv version of this paper \cite{Arxiv}.

\smallskip
\noindent
{\bf Acknowledgement} The first author gratefully acknowledges support by the Institut Mittag-Leffler (Djursholm,
Sweden) during the program Graphs, Hypergraphs and Computing.


\begin{thebibliography}{99}

\bibitem{AR} L. Asimow and  B. Roth, The rigidity of graphs, {\rm Trans. Amer. Math. Soc.} 245 (1978), 279-289.


\bibitem{BJ} A. R. Berg and T. Jord\'an, Algorithms for graph rigidity and scene analysis,
in {\rm Algorithms-ESA 2003}, Lecture notes in Computer Science
2832, 2003, 78-89.

\bibitem{BF} W. Bouma, I. Fudos, C. M. Hoffmann, J. Cai, and R. Paige, A Geometric
Constraint Solver, Computer Aided Design, 27(1995) 487-501.

\bibitem{B} T. Brylawski, Constructions, in {\rm Theory of Matroids},
ed. N. White, CUP, London, 1986, 127-223.

\bibitem{C} R. Connelly, Generic global rigidity,
{\rm Discrete Comput. Geom.} 33 (2005) 549-563.

\bibitem{Dunst} F. D. J. Dunstan, Matroids and submodular functions, {\rm Quart. J. Math.}, 27
(1976) 339-348.


\bibitem{Ematunion} J. Edmonds, Minimum partition of a matroid into independent subsets,
{\rm J. Res. Nat. Bur. Standards}, 69B (1965) 67-72.

\bibitem{E} J. Edmonds, Submodular functions, matroids, and certain polyhedra, in {\rm Combinatorial Structures and their Applications},
eds. R. Guy, H. Hanani, N. Sauer, and J. Sch\"{o}nheim, Gordon and Breach, New York, 1970, 69-87.

\bibitem{EF} J. Edmonds and D. R. Fulkerson, Transversals and matroid partition,
{\rm J. Res. Nat. Bur. Standards}, 69B (1965) 147-153.


\bibitem{F} {A. Frank},
{\rm Connections in combinatorial optimization}, Oxford Lecture
Series in Mathematics and its Applications, 38, Oxford University
Press, Oxford 2011.

\bibitem{FG} A. Frank and A. Gy\'arf\'as, How to orient the edges of a graph?, in {\rm Combinatorics},
(Keszthely), Coll. Math. Soc. J. Bolyai 18, North Holland 1976, 353-364.

\bibitem{GW} H. N. Gabow and H. H. Westermann, Forests, frames, and games: algorithms for matroid sums and applications, {\rm Algorithmica}, 7
(1992) 465-497.

\bibitem{GC} X.-S.Gao and S.-C. Chou, Solving geometric constraint systems. II A symbolic approach and
decision of rc-constructability, Computer-Aided Design 30 (1998), 115-122.

\bibitem{GH} S. D. Guest and J. W. Hutchinson, On the determinancy of repetitive structures,
{\rm J. Mechanics and Physics of Solids}, 51 (2003) 383-391.

\bibitem{I} H. Imai, On combinatorial structures of line drawings of polyhedra, {\rm Disc. Appl. Math.} 10 (1985) 79-92.



\bibitem{JH} D. J. Jacobs and B. Hendrickson, An algorithm for two-dimensional rigidity percolation:
the pebble game, {\rm J.Comput.Phys.} 137 (1997) 346-365.

\bibitem{JJ} B. Jackson and T Jord\'an, Connected rigidity matroids and unique
realizations of graphs, {\rm J. Combinatorial Theory(B)}, 94,
(2005), 1-29.

\bibitem{Arxiv} B. Jackson and J. Owen, A characterisation of the generic rigidity of 2-dimensional point-line frameworks, http://arxiv.org/abs/1407.4675

\bibitem{KT} N. Katoh and S. Tanigawa, A rooted-forest partition
with uniform vertex demand, J. Comb. Optim. 24 (2012) 67-98.


\bibitem{L} G. Laman, On graphs and rigidity of plane skeletal structures,
{\rm J. Engineering Math.} 4 (1970), 331-340.

\bibitem{LS} A. Lee and I. Streinu, Pebble game algorithms and sparse graphs,
{\rm Discrete Mathematics} 308 (2008) 1425-1437.

\bibitem{LY} L. Lov\'asz and Y. Yemini, On generic rigidity in the plane,
{SIAM J. Algebraic Discrete Methods} 3 (1982)  91-98.



\bibitem{O} J. C. Owen, Algebraic solution for geometry from dimensional constraints,
{\rm ACM Symposium on Foundations in Solid Modeling} (1991),
397-407.

\bibitem{OO} J. C. Owen, Constraints on simple geometry in two and three dimensions,
{\rm J. Comput. Geom. Appl.} 6 (1996) 421.

\bibitem{PP} J. Pym and H. Perfect, Submodular functions and independence
structures, J. Math. Anal. Appl. 30  (1970) 1-31.

\bibitem{S} K. Sugihara, Detection of structural inconsistency in systems of equations with degrees of freedom and its applications, {\rm Disc. Appl. Math.} 10 (1985) 312-328.



\bibitem{W} W. Whiteley, The union of matroids and the rigidity of frameworks,
{\rm SIAM J. Disc. Math.} 1 (1988) 237-255.


\bibitem{Wscene} W. Whiteley, A matroid on hypergraphs with applications in scene analysis and geometry,
{\rm Disc. Comput. Geom.} 4 (1989) 75-95.


\bibitem{WW} W. Whiteley, Matroids and rigid structures, in
{Matroid Applications} ed. Neil White, Encyclopedia of Mathematics
and its applications 40 (1992) 1-51.

\bibitem{WO} W. Whiteley, Representing Geometric Configurations: Computational Approaches,
in {\rm Learning and Geometry} eds. D. W. Kuecher and C. H. Smith,
Birkhauser 1996 143-177.

\bibitem{Wchapter} {W. Whiteley}, {Some matroids from discrete applied
geometry}, in {\rm Matroid Theory}, J. E. Bonin, J. G. Oxley, and B.
Servatius, Eds. American Mathematical Society, Contemporary
Mathematics 197 1996 171-313.


\bibitem{Y} L. Yang, Distance Coordinates Used in Geometric Constraint Solving,
in {\rm Automated Deduction in Geometry}, Lecture Notes in Computer
Science 2930  2004  216-229.

\end{thebibliography}
\end{document}